\documentclass[12pt,draftcls,onecolumn]{IEEEtran}

\IEEEoverridecommandlockouts                              

\makeatletter
\let\NAT@parse\undefined
\makeatother

\usepackage{cite}
\usepackage{amsmath,amssymb,amsfonts,amsthm}
\usepackage{algorithmic}
\usepackage{graphicx}
\usepackage{algorithm,algorithmic}
\usepackage{hyperref}
\hypersetup{hidelinks=true}
\usepackage{textcomp}
\usepackage{easyReview}
\usepackage{caption}
\usepackage{subcaption}
\usepackage{bbm}
\usepackage{comment}
\usepackage[export]{adjustbox}

\newtheorem{Def}{Definition}
\newtheorem{prop}{Proposition}
\newtheorem{thm}{Theorem}

\newtheorem{cor}{Corollary}
\newtheorem{rem}{Remark}
\newtheorem{ex}{Example}

\newcommand{\cb}[1]{{\color{blue}{#1}}}
\newcommand{\cred}[1]{{\color{red}{#1}}}

\def \b {\beta}
\def \O {\Omega}

\def \m {\mathfrak{m}}
\def \<{\langle}
\def \>{\rangle}

\def \w {\omega}
\def \t {\theta}
\def \a {\alpha}
\def \s {\sigma}

\def \r {\rho}

\def \O {\Omega}
\def \S {\Sigma}

\def \cP {\mathcal{P}}
\def \cF {\mathcal{F}}

\def \cM {\mathcal{M}}

\def \cK {\mathcal{K}}

\title{\LARGE \bf

Distributional Control of Ensemble Systems

}

\author{Jr-Shin~Li and Wei Zhang
\thanks{*This work was supported by the Air Force Office of Scientific Research under the award FA9550-21-1-0335.}
\thanks{J.-S. Li is with the Department of Electrical and Systems Engineering and the Division of Computational \& Data Sciences, Washington University in St. Louis, St. Louis MO, 63130 USA (e-mail:  jsli@wustl.edu).}%
\thanks{W. Zhang is with the Department of Electrical and Systems Engineering, Washington University in St. Louis, St. Louis MO, 63130 USA (e-mail: wei.zhang@wustl.edu).}%
}

\begin{document}

\maketitle

\begin{abstract}
    Ensemble control offers rich and diverse opportunities in mathematical systems theory. In this paper, we present a new paradigm of ensemble control, referred to as distributional control, for ensemble systems. We shift the focus from controlling the states of ensemble systems to controlling the output measures induced by their aggregated measurements. To facilitate systems-theoretic analysis of these newly formulated distributional control challenges, we establish a dynamic moment kernelization approach, through which we derive the distributional system and its corresponding moment system for an ensemble system. We further explore optimal distributional control by integrating optimal transport concepts and techniques with the moment representations, creating a systematic computational distributional control framework. 
\end{abstract}


\section{Introduction}
Ensemble systems, consisting of large populations of structurally similar dynamic units with heterogeneous dynamics, are prevalent in nature, societies, and engineered infrastructures. The ability to finely control the collective behavior within such large-scale systems is recognized as a fundamental and pivotal step, enabling diverse applications across various domains from quantum mechanics and neuroscience to robotics \cite{Glaser98,Li_PNAS11,ching_13_control,Kiss2007,Sun2020}. These \emph{ensemble control} tasks, however, present a two-fold challenge, as control and observation can only be conducted at the population level. These limitations consequently give rise to new and unconventional problems that extend beyond the capabilities of canonical tools in modern systems theory.

In recent years, substantial research efforts have been dedicated to addressing ensemble control problems. These efforts have primarily focused on investigating fundamental properties such as ensemble controllability, reachability, observability, and synchronizability in ensembles of isolated or networked systems \cite{Li_TAC09,Belhadj15,Chen_Automatica19,Boscain_SICON18,Boscain_JDE22,Dirr2016,Dirr2021,Schonlein2021,Schonlein2022,Zeng_TAC15,Chen_Automatica20_Bloch,Kiss2002}.
These investigations have led to the development of novel theoretical approaches that integrate tools from both closely related and remote areas. For example, concepts and techniques from polynomial approximation, spectral theory, complex analysis, graph theory, and Lie algebra have been bridged to study ensemble controllability for linear, bilinear, and some forms of nonlinear ensemble systems \cite{Li_PRA_2006,Li_TAC11,Li_TAC13,Helmke_SCL14,Zeng_SCL16,Chen_Automatica19,Chen_MCSS19,Li_SICON21}. In addition, ideas and tools from tomography, probability theory, and kernel methods have been synergistically devised and applied to the study of ensemble observability for linear ensemble systems \cite{Zeng_TAC15,Zeng_TAC16,Zeng_Automatica19,Zeng_CDC16,Chen_Automatica20_Bloch,Miao_thesis}. Recently, the ensemble system formulation has been extended to study partial differential equation (PDE) systems for the design of stabilizing boundary control laws \cite{Krstic_TAC24}. 
Of parallel importance, theory-driven and learning-based numerical methods have been extensively developed to enable effective and efficient design and computation of feasible and optimal ensemble control inputs \cite{GRAPE,Li_PNAS11,Li_JCP11,Dong_PRA14,Li_CDC15_Fourier,Gong_SICON16,Brent_JCP06,Li_ACC20,Li_ACC2024,Zlotnik_ACC24}.

Notably, these successful developments were achieved primarily through the analysis of state-space models that describe the time evolution of ensemble systems without exploiting their measurements. This approach is partly necessitated by the underactuated nature of ensemble systems, where obtaining state feedback information for each individual system in a large ensemble is either arduous or impossible.

In this paper,  we introduce the concept and formulation of the \emph{distributional control} problem arising in ensemble control systems. This new paradigm in ensemble control is motivated by the availability of population-level measurements, which we refer to as \emph{aggregated measurements}, and is concerned with controlling output distributions induced by these aggregated measurements. To put our study into a formal setting, we introduce the output measure associated with an ensemble system and focus on controlling the dynamics of these output measures. To facilitate systems-theoretic analysis of distributional control problems, we propose a dynamic moment kernelization approach that generates moment representations of the output measures of an ensemble system. Through this transformation, we derive the distributional system defined on the space of output measures and its corresponding moment system. We further study optimal distributional control by leveraging the concepts and techniques of optimal transport (OT), combined with the developed moment representations, to synthesize an OT-enabled distributional control method. 

It is worth noting that distributional control in this context is distinct from classical stochastic control or control of distributions. 
In these areas, the focus is on controlling a system driven by stochastic processes, such as Poisson counters or Winner processes \cite{Liberzon_SICON20,Brockett2000}, or on controlling a population of identical particles whose initial condition follows a specific probability distribution \cite{Brockett2012,Zeng_CDC16,Chen_TAC16_I,Chen_TAC16_II,Chen_TAC16_OT}. 


This paper is organized as follows. In Section \ref{sec:DC}, we introduce the notion of distributional control of the output distributions of ensemble systems using a measure-theoretic formulation and investigate the relationship between distributional and ensemble control in terms of the controllability properties of ensemble systems. Section \ref{sec:DMK}, which contains the main theoretical development, is devoted to the dynamic moment kernelization approach to distributional control. Here, we introduce moments of ensemble output measures and derive the moment system corresponding to a distributional system. 
In Section \ref{sec:OT_EC}, we study optimal distributional control through a systematic integration of OT techniques and demonstrate the performance of this method using illustrative and practical examples.

\section{Distributional Control of Ensemble Systems}
\label{sec:DC} 
In this section, we introduce a novel concept and a systems-theoretic formulation of a distinct class of problems in ensemble control, referred to as distributional control. 
These problems arise from the control of output distributions (measures) that are induced by population-level measurements of ensemble systems. We begin with the state-space description of ensemble control systems defined on a function space. Then, we shift our focus to the case where aggregated measurements of ensemble systems are available, which motivates the distributional control problem.

\subsection{Ensemble Control Systems Defined on Function Spaces}
\label{sec:controllability}
An \emph{ensemble system} refers to a sizable population of \emph{heterogeneous} dynamic units that share the same structure but exhibit distinct dynamics governed by different values of system parameters \cite{Li_TAC09}. In practice, an ensemble system may consist of a finite or an infinite number of systems, e.g., a continuum at the limit. A fundamental challenge in addressing ensemble systems lies in their inherent limitations, where control and observation can only be achieved at the population level. More specifically, ensemble systems are underactuated, as all dynamic units in an ensemble often receive a common source of control inputs, which appears as a restriction in many applications \cite{Li_PNAS11,Becker_TRO12}. Additionally, only aggregated measurements, such as data snapshots or images, are available, making state feedback for each individual system unattainable due to limitations in sensing capability \cite{Li_PRA_2006,Narayanan_ACC19,Li_ACC20_Learning,Dong_PRA22}. 

Mathematically, an ensemble system can be formulated as a parameterized control system defined on a function space, 
given by 
\begin{equation}
    \label{eq:ensemble_system}
    \frac{d}{dt}x(t,\beta) = F(x(t,\beta), \beta, u(t)).
\end{equation}
This parameterized system is indexed by the parameter $\b\in\Omega\subset\mathbb{R}^d$ varying on a compact subset $\O$. Each individual system evolves on a manifold $M$, i.e., $x(t,\b)\in M$ for all $t\in\mathbb{R}$ and $\b\in\O$, controlled by the common input $u\in U$ defined by $u:[0,T]\to\mathbb{R}^p$, where $U$ is a set of measurable functions. Therefore, at each time  $t\in\mathbb{R}$, the state $x(t,\cdot)$ is an $M$-valued function defined on $\O$, and the state-space of this ensemble system is a space of $M$-valued functions defined on $\O$, denoted by $\cF(\O,M)$. Analogous to classical control systems, steering an ensemble system from an initial state $x_0\in\mathcal{F}(\Omega,M)$ to a desired target state $x_F\in\mathcal{F}(\Omega,M)$ at a prescribed time $T>0$ is of fundamental significance and practical relevance \cite{Glaser98,Li_PNAS11,ching_13_control,Li_NatureComm16}. The primary obstacle in achieving this transfer lies in the restriction of utilizing a parameter-independent control input $u(t)\in\mathbb{R}^p$, i.e., a broadcast open-loop law. The ability to accomplish such a desired transfer between any pair $(x_0,x_F)$ is quantified by the notion of ensemble controllability.

 \begin{Def}[Ensemble controllability]
    \label{def:ensmeble_controllability}
    The system in \eqref{eq:ensemble_system} is said to be \emph{ensemble controllable} on $\mathcal{F}(\Omega,M)$ if, for any $\varepsilon>0$ and any initial state $x_0\in\mathcal{F}(\Omega,M)$, there exists a measurable control function $u(t)$ that steers the ensemble from $x_0$ into the $\varepsilon$-neighborhood of a desired target state $x_F\in\mathcal{F}(\Omega,M)$ at a finite time $T>0$. This implies $d_F(x(T,\cdot),x_F(\cdot))<\varepsilon$, where $d_F:\mathcal{F}(\Omega,M)\times\mathcal{F}(\Omega,M)\rightarrow\mathbb{R}$ is a metric on $\mathcal{F}(\Omega,M)$, and the final time $T$ may depend on $\varepsilon$.
\end{Def}

This fundamental property is defined in relation to the chosen topology, so that the underlying metric defined through $d_F$ plays a crucial role in the analysis of ensemble controllability \cite{Li_SICON21}. It is worth noting that ensemble controllability is a notion of approximate controllability \cite{Helmke_SCL14}.

\subsection{Ensemble Systems with Aggregated Measurements}
\label{sec:controllability}
\begin{figure}[h]
    \begin{center}
    \includegraphics[width=0.8\linewidth,keepaspectratio]{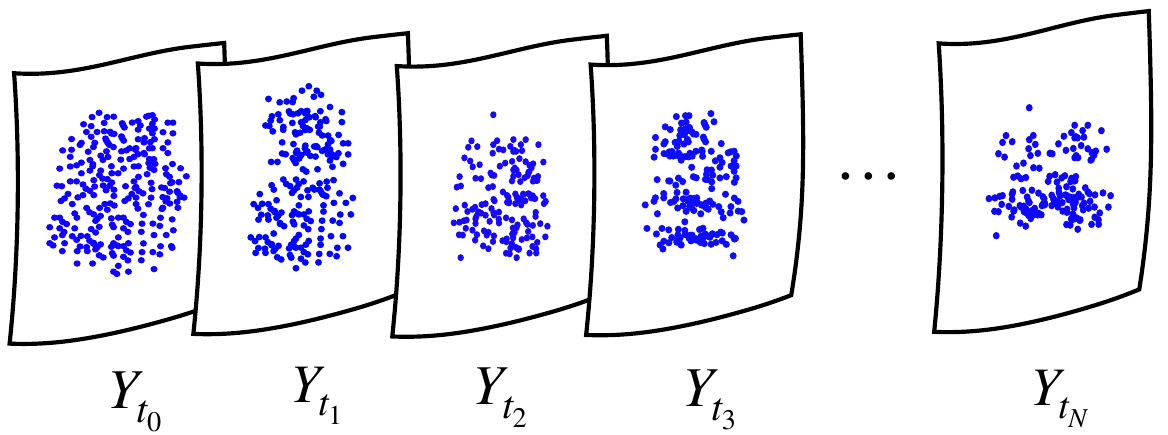}
    \end{center}
    \caption{\small Sample aggregated measurements $Y_{t_k}$ at time $t_k$ for $k=0,1,\ldots,N$ of an ensemble system.} 
    \label{fig:AM}
    \end{figure}

In practice, observations for ensemble systems are limited at the population level. The collected measurement data are referred to as \emph{aggregated measurements}, as illustrated in Figure \ref{fig:AM}. Together with the available aggregated measurements, an ensemble control system with dynamics described by \eqref{eq:ensemble_system} can be modeled as
\begin{align}
    \label{eq:ECS}
    \Sigma:
    \begin{cases}
        \frac{d}{dt}x(t,\b) = F(x(t,\b),\b,u(t)), \\ 
        \qquad \ \ \, Y_t = h\circ x_t(\O).
    \end{cases}
\end{align}
In this model, the aggregated measurements $Y_t$ are contained in the \emph{output space} $N$, which is assumed to be a Polish metric space, i.e., a separable and complete metric space \cite{Bogachev07}. 
The map $h:\mathcal{B}(M)\rightarrow\mathcal{B}(N)$ is the output function observing the states of all the systems in the ensemble at the population level as a set $Y_t=h(x_t(\O))$ at each time $t$, 
where $x_t(\O)=\{x_t(\beta)\in M:\beta\in\Omega\}$, and $\mathcal{B}(M)$ and $\mathcal{B}(N)$ denote the Borel $\sigma$-algebras on $M$ and $N$, respectively. 

An emergent critical observation from this model is that, despite the deterministic nature of this ensemble system, the aggregated measurement $Y_t$ reveals distributional information inherent in the system dynamics. Specifically, $Y_t$ induces a time-dependent measure, e.g., a probability distribution, defined on the observation space $N$. This interpretation gives rise to the concept of \emph{distributional control}, which pertains to the regulation of dynamic patterns within an ensemble system, e.g., from $Y_{t_0}$ to $Y_{t_N}$ as depicted in Figure \ref{fig:AM}, and plays a pivotal role in numerous emerging applications such as synchronization engineering and quantum information sciences \cite{Strogatz2000,Silver85}.

\begin{ex}[Probability distributions induced by aggregated measurements] \label{ex:distribution}
    \rm Consider the scalar linear ensemble system defined on $\mathcal{F}(\Omega,\mathbb{R})$, given by,
    \begin{align}
    \label{eq:Sigma1}
    \Sigma_1:
        \begin{cases}
            \frac{d}{dt}x(t,\beta)=\beta x(t,\beta)+u(t), \\ 
            \qquad\ \ \, Y_t =\mathbf{1}_{[0,\infty)}\circ x_t(\O), 
        \end{cases}
    \end{align}
    where the output function $h=\mathbf{1}_{[0,\infty)}$ is the indicator function of $[0,\infty)$. The system parameter $\b$ takes values in $\Omega\subseteq[0,1]$, and $h:\mathbb{R}\rightarrow\{0,1\}$ generates binary outputs given by $h(x_t(\beta))=0$ if $x_t(\beta)<0$, and $h(x_t(\beta))=1$ otherwise. As a result, the aggregated measurements only contain the binary values, i.e., $Y_t=\{0,1\}$, for all $t$. Then, the frequency of occurrences of 0 or 1 in each $Y_t$ will induce a time-varying probability distribution. 
    
    To embody this new ``distributional'' viewpoint with a tangible example, we select the initial condition $x_0(\beta)=\beta-a$ for some constant $a\in[0,1]$ and choose the control input $u(t)=0$ for all $t$. Then, the solution of the system in \eqref{eq:Sigma1} is $x_t(\beta)=e^{t\beta}x_0(\beta)$, which yields 
    \begin{align*}
       \begin{cases}  x_t(\beta)<0,\ \text{if } \beta\in[0,a), \\  x_t(\beta)\geq0,\ \text{if } \beta\in[a,1],
        \end{cases}\text{and}\ \  
        y_t(\beta)
        =\begin{cases} 0,\ \text{if }\beta\in[0,a), \\ 1,\ \text{if }\beta\in[a,1],
        \end{cases}
    \end{align*}
    where $y_t(\b)\doteq h(x_t(\b))$.
    It becomes evident that $h$ induces a probability distribution of $Y_t=\{y_t(\b):\b\in\O\}$, which is defined by the probability measure $\mu_t$ on $\{0,1\}$ satisfying $\mu_t(\{0\})=a$ and $\mu_t(\{1\})=1-a$. In fact, $\mu_t$ is a Bernoulli distribution, i.e., $\mu_t\sim\rm{Bernoulli}(a)$.

    Notice that the probabilities $a$ and $1-a$ are essentially the Lebesgue measure $\lambda$ (length) of the intervals $[0,a)$ and $[a,1]$ defined on $\O$. Specifically, they are the respective preimages of $0$ and $1$ under the function $y_t$. This allows us to define $\mu_t$ 
    in terms of the pushforward of $\lambda$ by $y_t$ as $\mu_t(\{0\})=\lambda([0,a))=\lambda(y_t^{-1}(\{0\}))=(y_t)_\#\lambda(\{0\})$ and $\mu_t(\{1\})=\lambda([a,1])=\lambda(y_t^{-1}(\{1\}))=(y_t)_\#\lambda(\{1\})$. 
    This pushforward representation can be readily generalized to the ensemble system $\S$ in \eqref{eq:ECS}.
\end{ex}

\begin{Def}[Output measure]
    \label{def:output_measure}
    Consider the ensemble system $\Sigma$ in \eqref{eq:ECS} and let $\lambda$ be a Borel measure on $\Omega$. The \emph{output measure} $\mu_t$ induced by the aggregated measurements $Y_t$ of $\Sigma$ is a Borel measure on $N$, defined by $\mu_t=(y_t)_\#\lambda$ with $y_t(\b)= h(x_t(\b))$.
\end{Def}

Note that $\mu_t$ satisfies $\mu_t(B)=\lambda\big(y_t^{-1}(B)\big)$ for any Borel set $B\subseteq N$, and 
this can be expressed as $\int_N fd\mu_t=\int_\Omega f\circ y_t d\lambda$ for any real-valued Borel measurable function $f$ on $N$ \cite{folland_13_real}. In particular, if $f$ is the constant function 1, then $\mu_t(N)=\int_N d\mu_t=\int_\Omega d\lambda=\lambda(\Omega)$. This implies that the total measure of $\mu_t$ is time-invariant and always equal to that of $\lambda$. Therefore, without loss of generality, we can assume that $\mu_t$ is a probability measure for all $t$. A natural choice of $\lambda$ leading to this assumption is a normalized Riemannian volume measure; that is, $d\lambda=\frac{\omega_g}{{\rm vol}_g(\Omega)}$, with $\omega_g$ being the volume form with respect to a Riemannian metric $g$ on $\O$. The compactness of $\Omega$ guarantees ${{\rm vol}_g(\Omega)}=\int_\Omega\omega_g<\infty$. More specifically, given a coordinate chart $(\beta_1,\dots,\beta_d)$ on $\Omega$, the metric tensor $g$ is determined by its local matrix representation $g=(g_{ij})\in\mathbb{R}^{d\times d}$, yielding $\omega_g=\sqrt{\det(g)}d\beta_1\cdots d\beta_d$ \cite{villani2009}. 

\subsection{Pattern Control in Ensemble Systems}
\label{sec:pattern_controllability}
Beyond controlling the time-evolution of an ensemble system, which has been the main focus of previous works in the field of ensemble control, the manipulation of time-varying output measures induced by aggregated measurements introduces a novel and rich class of distributional control problems that have not been explored in the literature. This emerging perspective is particularly pertinent to a broad range of applications, where dynamic pattern formation or distribution shaping is of fundamental and practical importance, such as synchronization in complex networks, targeted coordination in robot swarms, and pattern regulation in spiking neurons \cite{Strogatz2000,Becker_TRO12,Li_ACC20_Learning}.
The ability to achieve this control task can be quantified by a new concept of pattern controllability.


\begin{Def}[Pattern controllability]
    \label{def:pattern_controllability}
    Given the ensemble system $\Sigma$ in \eqref{eq:ECS} and the output measure $\mu_t\in\mathcal{P}(N)$ induced by its aggregated measurements $Y_t$, 
    $\Sigma$ is said to be \emph{pattern controllable} on $\mathcal{P}(N)$ if for any initial output measure $\mu_0\in\mathcal{P}(N)$ at time $t=0$ and any $\varepsilon>0$, there exists a measurable control function $u(t)$ that steers $\mu_0$ to be within the $\varepsilon$-neighborhood of a desired output measure $\mu_F$ at a finite time $T>0$, i.e., $d_P\big(\mu_T,\mu_F\big)<\varepsilon$, where $\cP(N)$ is the space of output measures, and $d_P:\mathcal{P}(N)\times\mathcal{P}(N)\rightarrow\mathbb{R}$ is a metric on $\mathcal{P}(N)$. 
\end{Def}

\begin{rem}
    Recall that the output space $N$ of the ensemble $\Sigma$ in \eqref{eq:ECS} is a Polish metric space. Thus, $\mathcal{P}(N)$ with the weak topology is metrizable, and every element $\mu\in\mathcal{P}(N)$ can be represented in the form of $\mu=f_\#\lambda$ for some $f:\Omega\rightarrow N$ \cite{Bogachev07}. 
\end{rem}


\begin{thm}
\label{thm:ensemble_pattern_controllability}
    If the ensemble system $\Sigma$ in \eqref{eq:ECS} is ensemble controllable on $\mathcal{F}(\Omega,M)$ and the output function $h:M\rightarrow N$ is continuous and surjective, then $\Sigma$ is pattern controllable on $\mathcal{P}\big(N\big)$.
\end{thm}
\begin{proof}
    Let $\mu_0\in\mathcal{P}(N)$ and $\mu_F\in\mathcal{P}(N)$ be a pair of initial and target output measures. Then, for any $\mu\in\mathcal{P}(N)$ and a given $\varepsilon>0$ such that $d_P(\mu_F,\mu)<\varepsilon$, there exist $\nu_F,\nu\in\mathcal{P}(M)$ and a $\delta>0$ satisfying $h_\#\nu_F=\mu_F$, $h_\#\nu=\mu$, and $d_P(\nu_F,\nu)<\delta$. This is due to the surjectivity and continuity of $h_\#:\mathcal{P}(M)\rightarrow\mathcal{P}(N)$ following the same properties of $h:M\rightarrow N$ \cite{Bogachev07}. It remains to show the existence of a control input $u(t)$ that steers the system $\Sigma$ to $x_T\in\mathcal{F}(\Omega,M)$ satisfying $d_P(\nu_F,(x_T)_\#\lambda)<\delta$ in a finite time $T$. 


    To this end, using Skorokhod's representation theorem \cite{Bogachev07}, there exists an $x_F\in\mathcal{F}(\Omega,M)$ such that $\nu_F=(x_F)_\#\lambda$. This yields that a sequence of ensemble states, $x_{t_n}$, converging to $x_F$ 
    in distribution if and only if $(x_{t_n})_\#\lambda\rightarrow\nu_F$ weakly. Now, we define a metric $d_F$ on $\mathcal{F}(\Omega,M)$ by $d_F(f_1,f_2)=\int_\Omega\frac{\rho(f_1,f_2)}{1+\rho(f_1,f_2)}d\lambda$ for $f_1,f_2\in\mathcal{F}(\Omega,M)$, where $\rho$ is the distance on $M$ induced by a Riemannian metric. 
    Then, the convergence of $x_{t_n}$ in $d_F$, i.e., $d_F(x_{t_n},x_F)\rightarrow0$, implies its convergence in probability in probability and hence in distribution \cite{Billingsley95}. Equivalently, there exists a $\gamma>0$ such that $d_P((x_{t_n})_\#\lambda,\nu_F)<\delta$ whenever $d_F(x_{t_n},x_F)<\gamma$. As a result of ensemble controllability, there is a control input $u(t)$ steering $\Sigma$ from any initial condition into an $\gamma$-neighborhood of $x_F$, which in turn establishes the pattern controllability of $\Sigma$. 
\end{proof}

Theorem \ref{thm:ensemble_pattern_controllability} indicates that ensemble controllability is a sufficient condition for pattern controllability.  However, in general, this is not a necessary condition since distinct ensemble states can generate the same output measure. This can be directly observed through Example \ref{ex:distribution}, where the state functions $x_{t_1}\neq x_{t_2}$ for $t_1\neq t_2$, while the output measure $\mu_t$ remains the Bernoulli distribution, regardless of $t$. 

\section{Dynamic Moment Kernelization for \\ Ensemble Systems}
\label{sec:DMK}
The presented distributional control arising from ensemble systems forms a new class of problems in control theory, defined on a space of measures. These problems are distinct from conventional stochastic control problems, such as population control of identical particles governed by Liouville's or Schr\"{o}dinger equations \cite{Brockett2012,brockett_2000_stochastic}, and the control of probability densities described by Fokker-Planck equations for a stochastic system driven by noise processes \cite{Liberzon_SICON20}. To facilitate formal analysis and establish robust control design principles and methodologies to address this emerging field, we introduce the method of moment kernelization, which induces appropriate coordinate systems in which ensemble systems can be expressed in terms of succinct kernel representations.

\subsection{Moments of Ensemble Output Measures}
\label{sec:moment_ES}
As our focus in distributional control is on regulating the output measures associated with an ensemble system, it is advantageous to represent these measures in appropriate coordinates. To construct a specific coordinate system, we introduce the moment kernel transform. This transform assigns each output measure an infinite sequence of ``moments,'' yielding moment coordinates in the output measure space. 

To put this development into a rigorous mathematical setting, we further suppose that the output space $N$ is locally compact. Then, each output measure $\mu_t\in\mathcal{P}(N)$ is a Radon measure and hence admits a coordinate representation in terms of an infinite sequence (see Appendix\ref{appd:measure}). This can be constructed using a primal-dual pairing of the form,
\begin{equation}
    \label{eq:moment_transform}
    m_k(t)=\<\psi_k,\mu_t\>,
\end{equation}
where $\{\psi_k\}_{k\in\mathbb{N}}$ is a basis of $C_0(N)$, the space of continuous real-valued functions on $N$ vanishing at infinity. We refer to $m_k(t)$ as the \emph{$k^{\rm th}$ moment} and $m(t)=(m_0(t),m_1(t),m_2(t),\cdots)'$ as the \emph{moment sequence} of the output measure $\mu_t$.

\begin{rem}
\rm
\label{rem:duality}
    By defining the moments through the primal-dual pairing $\<\cdot,\cdot\>:C_0(N)\times\mathcal{P}(N)\rightarrow\mathbb{R}$, we treat an output measure $\mu_t\in\mathcal{P}(N)$ as a continuous linear functional on $C_0(N)$. The Riesz–Markov–Kakutani representation theorem then implies that the pairing is given by the integration \cite{folland_13_real},
    \begin{align}
        \label{eq:moment_integral}
        \<\varphi,\mu_t\>=\int_N\varphi d\mu_t.
    \end{align}
    Because $N$ is a locally compact Hausdorff space, 
    the primal-dual pairing between $C_0(N)$ and $\mathcal{P}(N)$ is well-defined (see Appendix\ref{appd:measure}).
    
\end{rem}

The definition in \eqref{eq:moment_transform} not only defines the moment coordinates, $m_k(t)$, of an output measure, $\mu_t$, but also identifies a fundamental relationship between the output measures of an ensemble system and their corresponding moment sequences.

\begin{thm}
\label{thm:one-to-one}
    Given the ensemble system $\Sigma$ in \eqref{eq:ECS}, the output measure $\mu_t$ is in one-to-one correspondence with its associated moment sequence $m(t)$.
\end{thm}
\begin{proof} 
    By the Riesz–Markov–Kakutani representation theorem, the integral operator $\psi\mapsto\int_N\psi d\mu_t$ is continuous on $ C_0(N)$ \cite{folland_13_real}. Hence, the moments $m_k(t)$ are well-defined for all $k\in\mathbb{N}$, as $|m_k(t)|=\big|\int_N\psi_kd\mu_t\big|\leq\sup_{y\in N}|\psi_k(y)|\int_Nd\mu_t=\sup_{y\in N}|\psi_k(y)|<\infty$. It remains to show that $\mu_t$ is uniquely determined by $m(t)$ and vice versa. This directly follows from the facts that $\mu_{t_n}\rightarrow\mu_t$ weakly if and only if $m(t_n)\rightarrow m(t)$ componentwise (see Theorem \ref{thm:weak_convergence} in Appendix \ref{appd:measure}), and that both of these limits are unique. 
\end{proof}

Theorem \ref{thm:one-to-one} further implies that $(\mathcal{P}(N),\mathcal{K})$ defines a (global) coordinate chart on $\mathcal{P}(N)$, where $\mathcal{K}:\mathcal{P}(N)\rightarrow\mathcal{M}(N)$ is the \emph{moment transform} defined by $\mu_t\mapsto m(t)$, and $\mathcal{M}(N)$ denotes the space of moment sequences associated with the output measures in $\mathcal{P}(N)$. 


\begin{cor}
    The moment transform $\mathcal{K}:\mathcal{P}(N)\rightarrow\mathcal{M}(N)$, defined by $\mu_t\mapsto m(t)$, is a homeomorphism with respect to the weak and the product topology on $\mathcal{P}(N)$ and $\mathcal{M}(N)$, respectively.  
\end{cor}
\begin{proof}
    From Theorem \ref{thm:one-to-one}, we know that $\mathcal{K}$ is a bijective map and that $\mu_{t_n}\rightarrow\mu_t$ weakly if and only if $m_k(t_n)\rightarrow m_k(t)$ for all $k\in\mathbb{N}$. Hence, proving the continuity of $\mathcal{K}$ and $\mathcal{K}^{-1}$ reduces to verifying that both $\mathcal{P}(N)$ and $\mathcal{M}(N)$ satisfy the first axiom of countability. This follows immediately from the fact that $\mathcal{P}(N)$ and $\mathcal{M}(N)$ are metric spaces \cite{Bogachev07}.
\end{proof}

In widely encountered practical applications where $N$ is compact, i.e., measurement values are bounded within a specific range, the space $C_0(N)$ coincides with $C(N)$, the space of all continuous real-valued functions on $N$. For example, if $N=[a,b]$, then $C_0([a,b])=C([a,b])$. Because the set of polynomials is dense in $C([a,b])$, this allows us to define the monomial moments of $\mu_t\in\mathcal{P}(N)$ by $m_k(t)=\int_a^by^kd\mu_t(y)$ with the basis $\psi_k(y)=y^k$, $k\in\mathbb{N}$. The one-to-one correspondence between $\mu_t$ and $m(t)=\big(m_k(t)\big)_{k\in\mathbb{N}}$ in this case revives the Hausdorff moment problem \cite{hausdorff_23_momentprobleme}.

\subsection{Moment Dynamics of Time-Varying Output Measures}
The homeomorphism between the spaces of output measures and moment sequences allows us to use moment sequences as a coordinate system to represent an ensemble system and the time-evolution of its output distribution. In the following, we introduce and derive \emph{the distributional system} associated with the ensemble system $\Sigma$ in \eqref{eq:ECS}, which describes the dynamics of output measures.


\subsubsection{Distributional systems on the space of output measures}
\label{sec:distributional_systems}
To construct a dynamical system of output measures $\mu_t$ defined on $\mathcal{P}(N)$, we equip the ensemble system $\Sigma$ with appropriate regularity conditions. 





\begin{prop}
    \label{prop:distribution_dynamics}
    Consider the ensemble system $\Sigma$ in \eqref{eq:ECS} and suppose that the output function $h$ is Lipschitz continuous on $M$, then the output measure $\mu_t$ is weakly differentiable 
    on $[0,T]$ almost everywhere. 
    \end{prop}
\begin{proof}
    We introduce a smooth ``test function'' $\varphi\in C^\infty_c(N)$ on the space of compactly supported real-valued smooth functions on $N$. Since $x_t$ is locally Lipschitz continuous in $t$, and $h$ is Lipschitz continuous, $\varphi\circ y_t=\varphi\circ h\circ x_t$ is locally Lipschitz continuous in $t$. Hence, it is differentiable with respect to almost every $t\in[0,T]$ and all $\b\in\O$. Specifically, this gives the time-derivative, for almost every $t\in[0,T]$, 
    \begin{align}
        \frac{d}{dt}\<\varphi,\mu_t\>&=\frac{d}{dt}\int_N\varphi d\mu_t=\int_\Omega\frac{d}{dt}(\varphi\circ y_t)d\lambda, \nonumber\\
        &=\int_\Omega(\nabla\varphi\circ y_t)\cdot \frac{d}{dt}(h\circ x_t) d\lambda, \label{eq:distribution_dynamics} \\
        &=\int_\Omega(\nabla\varphi\circ y_t)\cdot (h_*F\circ y_t)d\lambda, \nonumber 
\end{align}
where we used the relation $\mu_t=(y_t)_\#\lambda$ defined in Definition \ref{def:output_measure}, and $h_*F$ denotes the pushforward of $F$ by $h$. Note that the interchange of integration and differentiation in the second equality follows from the dominant convergence theorem, since $\nabla\varphi$ is compactly supported. 
\end{proof}

This weak differentiability and the continuity property (see Appendix\ref{appd:continuity}) 
then guarantee that $\mu_t$ is necessarily a weak solution of a certain differential equation system of the form,


\begin{align}
    \label{eq:distribution_system}
    \frac{d}{dt}\mu_t=\mathfrak{F}(t,\mu_t,u(t)),
\end{align}
where $\mu_t\in\cP(N)$, $\mathfrak{F}$ is a vector field on $\mathcal{P}(N)$, and $u(t)$ is the same control input for $\Sigma$. We refer to this system as the \emph{distributional system} associated with $\S$. It follows that the dynamic equation
\begin{align}
    \label{eq:distribution_system_integral}
    \frac{d}{dt}\<\varphi,\mu_t\>=\<\varphi,\frac{d}{dt}\mu_t\>=\<\varphi,\mathfrak{F}(t,\mu_t,u(t))\>
\end{align}
holds for any test function $\varphi\in C_c^\infty(N)$. 

\begin{ex}
    To illustrate the derivation of a distributional system related to an ensemble system described above, we revisit the scalar ensemble system $\Sigma_1$ in \eqref{eq:Sigma1} in the absence of control input, i.e., $\frac{d}{dt}x(t,\beta)=\beta x(t,\b)$ with $\b\in\O=[0,1]$, and consider the identity output function
    $h(x_t(\beta))=x_t(\b)$. In this case, the output measure is a finite Borel measure on $\mathbb{R}$, given by $\mu_t=(x_t)_\#\lambda$, where $x_t(\beta)=e^{t\beta}x_0(\beta)\doteq\Phi_t(\beta)x_0(\beta)$, and $\lambda$ is the Lebesgue measure on $[0,1]$. To explicitly calculate $\mu_t$, it is essential to determine its \emph{distribution function}, defined by $F_t(y)=\mu_t\big((-\infty,y]\big)=\lambda\big(\{\beta\in[0,1]:e^{t\beta}x_0(\beta)\leq y\}\big)$ \cite{Billingsley95}. For the sake of illustration, we choose the initial distribution $\mu_0$ as the point mass at $a>0$; equivalently, $x_0(\b)=a$ for $\b\in[0,1]$ almost everywhere. Then, we have 
    \begin{align*}
        F_t(y)=
        \begin{cases}
        0,\qquad\quad\ \, y<a,\\
        \frac{1}{t}\log\frac{y}{a},\quad a\leq y\leq ae^t,\\
        1,\qquad\quad\ \, y>ae^t.
        \end{cases}
    \end{align*}
    This gives
    \begin{align*}
        \frac{\partial}{\partial t}F_t(y)=
        \begin{cases}
        -\frac{1}{t^2}\log\frac{y}{a},\quad a\leq y\leq ae^t,\\
        0,\qquad\qquad\ \ y<a \text{ or } y>ae^t,
        \end{cases}
    \end{align*}
    which yields the system governing the evolution of $F_t$ as
    \begin{align}
    \frac{\partial}{\partial t}F_t(y)=-yF_t(y)\frac{\partial}{\partial y}F_t(y)=-\frac{1}{2}y\frac{\partial}{\partial y}F_t^2(y). \label{eq:distribution_system_cmf}
    \end{align}
    Note that $F_t^2$ is also a distribution function, and we denote its associated measure by $\mu_t^2$. Using the primal-dual notation, we have $F_t(y)=\<\mathbf{1}_{(-\infty,y]},\mu_t\>$ and $F_t^2(y)=\<\mathbf{1}_{(-\infty,y]},\mu_t^2\>$, where $\mathbf{1}_{(-\infty,y]}$ denotes the indicator function on $(-\infty,y]$ with $\mathbf{1}_{(-\infty,y]}(z)=1$ if $z\in(-\infty,y]$ and $\mathbf{1}_{(-\infty,y]}(z)=0$ otherwise. This leads to the measure-theoretic representation of the distributional system in \eqref{eq:distribution_system_cmf} as
      \begin{align*}
    \frac{\partial}{\partial t}\mu_t=-\frac{1}{2}I\nabla\cdot\mu^2_t, 
    \end{align*}
    where $I$ denotes the identity function on $\mathbb{R}$.

\end{ex}

\begin{rem}
    If the ensemble $\S$ is homogeneous, consisting of identical dynamic units, then the vector field $F$ is independent of $\beta$. In this case,
    the integro-differential equation in \eqref{eq:distribution_dynamics} is reduced to
    \begin{align*}
    \frac{d}{dt}\<\varphi,\mu_t\>=\int_\Omega\nabla\varphi\cdot h_*Fd\mu_t=-\<\varphi,\nabla\cdot(\mu_th_*F)\>,
    \end{align*}
    where the last equality follows from the definition of the derivative of Schwartz distributions \cite{folland_13_real}. This implies that $\mu_t$ satisfies the \emph{continuity equation},
    \begin{align*}
        \frac{\partial}{\partial t}\mu_t+\nabla\cdot(\mu_th_*F)=0,
    \end{align*}
    which is widely used as a differential equation representation of conservation laws. This equation has broad applications ranging from fluid dynamics and quantum mechanics to optimal transport \cite{villani2009,Brockett2012} . 
\end{rem}

\subsubsection{Distributional systems in the moment coordinates} 
The distributional system derived in \eqref{eq:distribution_system} can be represented in terms of moment coordinates. This representation facilitates a transparent analysis of intricate distributional control problems.

Using the moment transform defined in \eqref{eq:moment_transform} and the dynamic equation derived in \eqref{eq:distribution_system_integral}, we obtain the moment dynamics obeying
\begin{align}
    \frac{d}{dt}m_k(t)&=\frac{d}{dt}\<\varphi_k,\mu_t\>=\<\psi_k,\mathfrak{F}(t,\mathcal{K}^{-1}m(t),u(t))\>\nonumber\\
    &\doteq \bar F_k(t,m(t),u(t)) 
    \label{eq:moment_system_k}
\end{align}
for all $k\in\mathbb{N}$. 
Note that here the test functions $\varphi_k$ are chosen to be the basis functions, $\psi_k$, of $C_0(N)$. This is feasible because $C_c^\infty(N)$ is dense in $C_0(N)$ under the supremum norm topology \cite{folland_13_real}.
%
Therefore, the moment system associated with the distributional system \eqref{eq:distribution_system} of $\S$ is then given by
\begin{align}
    \label{eq:moment_system}
    \frac{d}{dt}m(t) &=\bar F(t,m(t),u(t)) \\
    &= \big(\mathcal{K}_*\mathfrak{F}\big)(t,m(t),u(t)) = \mathcal{K}\big(\mathfrak{F}(t,\mathcal{K}^{-1}m(t),u(t))\big), \nonumber
\end{align}
where $m(t)=(m_0(t),m_1(t),m_2(t),\ldots,)'$, and $\bar F=\mathcal{K}_*\mathfrak{F}$ defined on the moment space $\mathcal{M}(N)$ is the pushforward of the vector field $\mathfrak{F}$ by the moment transform $\cK$. 
Specifically, $\bar F=\mathcal{K}\circ\mathfrak{F}\circ\mathcal{K}^{-1}$ is, in fact, a change of coordinates under the moment transformation $\cK$. 

\section{Optimal Distributional Control}
\label{sec:OT_EC}
The theoretical developments in the previous sections have laid a foundational framework for introducing and deriving the distributional system induced by time-varying aggregated measurements of an ensemble system. The next critical phase is to understand how control of distributional systems can be achieved. In this section, we will capitalize on the concepts and techniques of optimal transport to devise a systematic approach for controlling distributional systems.

\subsection{Time-Dependent Optimal Transport}
\label{sec:TDOT}
Optimal transport (OT) is concerned with transporting one probability measure to another at minimal cost \cite{villani2009}. This objective directly connects with the distributional control problem introduced for ensemble systems. We will leverage this observation to establish an OT-enabled distributional control paradigm that provides a principled approach to optimal distributional control and inspires a new perspective on the application of OT in studying ensemble systems.

In the classical setting of OT, the transportation between two probability measures, say $\r_0$ and $\r_1$ in $\mathcal{P}(N)$, is conceptualized as a one-step process. This context has recently been enriched with a dynamic interpretation through the introduction of displacement interpolation (DI) \cite{villani2009}. The main idea of DI is to interpret the OT trajectory from $\r_0$ to $\r_1$ as a time-dependent parameterization, denoted $\r_t$ for $t\in[0,1]$. Specifically, along this trajectory, the transport between any pair of distributions, $(\r_\tau,\r_\s)$, is optimal for $\tau,\sigma\in[0,1]$ with $\tau<\sigma$. In other words, the transport cost from $\r_\tau$ to $\r_\sigma$ is minimal among all possible transports. 

Mathematically, this interpolation procedure can be formulated as an infinite-dimensional constrained optimization problem over the space of flows on $N$, 
given by
\begin{align}
    J_{DI}=\min_{\{\Phi_t:N\rightarrow N\}_{0\leq t\leq1}}&\quad \int_N c(\Phi_t(y))d\r_0(y), \nonumber\\
    {\rm s.t.}\quad\quad\ & \quad \Phi_0=I,\quad (\Phi_1)_\#\r_0=\r_1,\label{eq:DI}
\end{align}
where $c:C([0,1],N)\rightarrow\mathbb{R}$ is the cost functional on the space of 
continuous curves on $N$, and $I:N\rightarrow N$ denotes the identity map on $N$. 
The solution to this optimization problem, i.e., the minimizer $\Phi_t^*$, depicts the OT trajectory from $\r_0$ to $\r_1$ by $\r_t^*=(\Phi_t^*)_\#\r_0$ \cite{villani2009}. 

In the case where $N=\mathbb{R}$ and $c(\gamma_t)=\int_0^1|\dot\gamma_t|^pdt$ for $p\geq1$ with $\dot\gamma_t$ being the time-derivative of the curve $\gamma_t$, the DI is given by 
\begin{align}
    \r^*_t=[(1-t)I+tG_1^{-1}\circ G_0]_\#\r, \label{eq:DI}
\end{align}
where $G_0$ and $G_1$ are the cumulative distribution functions of $\r_0$ and $\r_1$, respectively, and $G_1^{-1}(x)=\sup\{y\in\mathbb{R}:G_1(y)\leq x\}$ is the \emph{generalized inverse} of $G_1$ \cite{McCann97}. This is referred to as \emph{McCann's interpolation}, where the minimal transport cost $J_{DI}$ coincides with the \emph{Wasserstein metric} (OT distance), $W_p(\r_0,\r_1)=\big(\inf\big\{\int_{\mathbb{R}^2}|x-y|^pd\r(x,y):\r\in\Gamma(\r_0,\r_1)\big\}\big)^{\frac{1}{p}}$, between $\r_0$ and $\r_1$ in $\mathcal{P}(\mathbb{R})$, where $\Gamma(\r_0,\r_1)$ denotes the space of probability measures on $\mathbb{R}^2$ with the marginals $\r_0$ and $\r_1$ \cite{villani2009}.
\subsection{Optimal Distributional Control using Moment Representations}
\label{sec:tracking}
As the connection between OT and distributional control is revealed, we will employ the proposed moment kernelization method to facilitate the integration of OT principles into the distributional control framework. To elaborate on this idea, we consider the control of output measures of the ensemble system $\S$ from $\mu_0$ to $\mu_1$. As described in Section \ref{sec:TDOT}, the optimal transport trajectory $\mu_t^*$, $t\in [0,1]$, from $\mu_0$ to $\mu_1$ is characterized through DI. Following \eqref{eq:moment_integral}, the $k^{\rm th}$ moment coordinate of $\mu^*_t$ is given by 
\begin{align}
    m^*_k(t)=\<\psi_k,\mu^*_t\>=\int_N\psi_kd(\Phi_t^*)_\#\mu_0=\int_N\psi_k\circ\Phi_t^*d\mu_0.\label{eq:OT_moment}
\end{align}
On the other hand, the moment system representation of the distributional system in \eqref{eq:distribution_system} was derived in \eqref{eq:moment_system}. Therefore, to achieve the desired transport from $\mu_0$ to $\mu_1$, it is essential to control this moment system such that the moment trajectory $m(t)$ meets the DI trajectory $m^*(t)$. In this way, this originally challenging problem of controlling output measures over 
$\cP(N)$ is reduced to an optimal tracking problem with the OT moment trajectory $m^*(t)=(m^*_0(t),m^*_1(t),m^*_2(t),\ldots,)'$ as the reference trajectory, given by
\begin{align}
    \min_{u:[0,1]\rightarrow\mathbb{R}^r} & \  \int_0^1 d_{\mathcal{M}}(m^*(t),m(t))dt\nonumber\\
    {\rm s.t.}\quad & \ \frac{d}{dt}m(t)=\bar F(t,m(t),u(t)), \label{eq:OT_tracking}
\end{align}
where $d_\mathcal{M}$ is a metric on the space of moment sequences $\mathcal{M}(N)$ that metrizes the product topology. A canonical choice is $d_{\mathcal{M}}(m^*(t),m(t))=\sum_{k=0}^\infty 2^{-k}|m^*_k(t)-m_k(t)|$ \cite{folland_13_real}.

To tackle this infinite-dimensional tracking problem, 
we conduct a finite-dimensional approximation of the moment system. 
Formally, let $P_q:\mathcal{M}\rightarrow\widehat{\mathcal{M}}_q$ be the projection onto the space of order-$q$ truncated moment sequences $\widehat{\mathcal{M}}_q\subset\mathbb{R}^{q+1}\otimes\mathbb{R}^{n}$, defined by $\widehat m^q(t)=P_qm(t)=\big(m_0(t),\cdots,m_q(t)\big)'$, then the system governing the dynamics of $\widehat m^q(t)$ 
is given by
\begin{align}
    \frac{d}{dt}\widehat m^q(t)&=\frac{d}{dt}P_q m(t)=\nabla P_q\Big(\frac{d}{dt}m(t)\Big)\nonumber\\
    &=(P_q)_*\bar F(t,P_qm(t),u(t))\nonumber\\
    &\doteq\widehat F^q(t,\widehat m^q(t),u(t))
    \label{eq:truncated_moment_system},
\end{align}
which is an $n(q+1)$-dimensional system defined on $\widehat{\mathcal{M}}_q$. 
Similar to the moment system presented in \eqref{eq:moment_system}, $\widehat F^q$ can be regarded as a change of coordinates via $\widehat F^q=P_q\circ\bar F\circ P_q'$, where $P_q'$, the adjoint operator of $P_q$, is in fact the inclusion map $\iota_q:\widehat{\mathcal{M}}_q\hookrightarrow\mathcal{M}$.




\begin{thm} {\it (Convergence of truncated moment sequences):}
\label{prop:moment_convergnece}
    Consider the moment system in \eqref{eq:OT_tracking} with the initial condition $m(0)$ and its order-$q$ truncated moment system in \eqref{eq:truncated_moment_system} with $\widehat m^q(0)=P_qm(0)$, both controlled by the same input $u(t)$. Then, the trajectory $\widehat m^q(t)$ converges to the trajectory $m(t)$ uniformly on $0\leq t\leq 1$; that is, $\sup_{t\in[0,1]}d_\mathcal{M}(\widehat m^q(t),m(t))\rightarrow0$ as $q\rightarrow\infty$.
\end{thm}

\begin{proof}
    For all $k\leq q$, we have 
    \begin{align*}
        \frac{d}{dt}\big(\widehat m^q_k(t)-m_k(t)\big)&=\widehat F^q_k(t,\widehat m^q(t),u(t))-\bar F_k(t, m(t),u(t))\\
        &\leq|\bar F_k(t,\widehat m^q(t),u(t))-\bar F_k(t, m(t),u(t))|\\
        &\leq Ld_{\mathcal{M}}\big(\widehat m^q(t),m(t)\big), 
    \end{align*}
     where $L$ is the Lipchitz constant of $\bar F$. Since the same inequality also holds for $m_k(t)-\widehat m^q_k(t)$, this implies that 
    \begin{align}
        \frac{d}{dt}\big|\widehat m^q_k(t)-m_k(t)\big|\leq Ld_{\mathcal{M}}\big(\widehat m^q(t),m(t)\big). \label{eq:bound_difference}
    \end{align}
    On the other hand, when $k>q$, $\widehat m_k^q(t)=0$. This gives $\widehat F_k^q(t,\widehat m^q(t),u(t))=0$ for all $t$, and thus \eqref{eq:bound_difference} also holds true for $k>q$. Consequently, we obtain
        \begin{align*}
        &\frac{d}{dt}d_{\mathcal{M}}\big(\widehat m^q(t),m(t)\big) \\
        &=\frac{d}{dt}\sum_{k=0}^\infty2^{-k}\big|\widehat m^q_k(t)-m_k(t)\big|=\sum_{k=0}^\infty2^{-k}\frac{d}{dt}\big|\widehat m^q_k(t)-m_k(t)\big| \\ 
        &\leq\sum_{k=0}^\infty2^{-k}Ld_{\mathcal{M}}\big(\widehat m^q(t),m(t)\big)=2Ld_{\mathcal{M}}\big(\widehat m^q(t),m(t)\big),
    \end{align*}
    where the exchange of the order of differentiation and summation follows from the monotone convergence theorem \cite{folland_13_real}. Gr\"{o}nwall's inequality then implies $d_{\mathcal{M}}\big(\widehat m^q(t),m(t)\big)\leq e^{2L}\cdot d_{\mathcal{M}}\big(\widehat m^q(0),m(0)\big)$ for all $t\in[0,1]$. Since $\widehat m^q(0)$ consists of the first $q+1$ components of $m(0)$, it converges to $m(0)$ componentwise as $q\rightarrow\infty$. It follows that $d_{\mathcal{M}}\big(\widehat m^q(0),m(0)\big)\rightarrow0$ and hence $d_{\mathcal{M}}\big(\widehat m^q(t),m(t)\big)\to 0$ uniformly for all $t\in[0,1]$. 
\end{proof}

The derived moment convergence property allows us to address the OT-tracking control problem in \eqref{eq:OT_tracking} using truncated moment systems. This approximation leads to the finite-dimensional optimal tracking problem, given by
\begin{align}
    \min_{u:[0,1]\rightarrow\mathbb{R}^r}&\  \int_0^1 d_{\mathcal{M}}\big(P_qm^*(t),\widehat m^q(t)\big)dt\nonumber\\
    {\rm s.t.}\ \ \, & \frac{d}{dt}\widehat m^q(t)= \widehat{F}^q(t,\widehat m^q(t),u(t)), \label{eq:OT_tracking_truncated}
\end{align}
which can be solved using various existing methods, such as pseudospectral \cite{Gong_SICON16,Li_PNAS11} or iterative methods \cite{zeng_18_computation,Li_CDC24}.

\begin{cor}[Convergence in optimal tracking cost]
\label{cor:moment_convergence}
    Let $u_q$ be an optimal OT-tracking control for the problem in \eqref{eq:OT_tracking_truncated}, and $J_q(u_q)$ be the respective tracking cost. Then, $J_q(u_q)$ converges to the minimal tracking cost, denoted $J^*$, of the optimal control problem in \eqref{eq:OT_tracking}, i.e., $J_q(u_q)\to J^*$ as $q\rightarrow\infty$. 
\end{cor}
\begin{proof}
   The proof is based on constructing a double sequence $(J_p(u_q))_{p,q\in\mathbb{N}}$ that embeds $(J_q(u_q))_{q\in\mathbb{N}}$ as the diagonal subsequence. One can then demonstrate that both iterated limits, $\lim_{p\rightarrow\infty}\lim_{q\rightarrow\infty}$ $J_p(u_q)$ and $\lim_{q\rightarrow\infty}\lim_{p\rightarrow\infty}J_p(u_q)$, exist and are equal to $J^*$, implying $J_q(u_q)\rightarrow J^*$ as $q\rightarrow\infty$ 
   (see Appendix \ref{appd:proof_cor_2} for the detailed proof).
\end{proof}

\subsection{Examples and Simulations}
\label{sec:examples}
In this section, we present several examples to illustrate the developed concepts and theory of distributional control. These examples demonstrate the utilization of dynamic moment kernelization and the OT-enabled formulations and techniques for solving distributional control problems.

\begin{ex}[Functional and pattern control] 
\rm In this example, we illustrate two distinct aspects of distributional control problems arising from different types of aggregated measurements. Through this example, we provide a detailed workflow for utilizing the proposed moment method, combined with the OT-inspired technique, for distributional control of ensemble systems. 

Given a multi-input scalar-valued linear ensemble system, 
    \begin{align}
    \label{eq:ECS_linear}
    \Sigma_2:
        \begin{cases}
            \frac{d}{dt}x(t,\beta)=\beta x(t,\beta)+\sum_{i=1}^p\beta^{i-1} u_i(t), \\ 
            \qquad\ \ \, Y_t=h\circ x_t(\O), 
        \end{cases}
    \end{align}
where $x_t=x(t,\b)\in L^2(\O)$, and $\O= [0,1]$ is endowed with the Lebesgue measure $\lambda$, we consider the problem of controlling the output measure $\mu_t$ induced by the aggregated measurements $Y_t$ from $\mu_0$ to $\mu_1$.

Case I \emph{(Aggregated measurements with labels):} We first consider the case where the output function encodes the label $\b$ of each individual system in the ensemble. Specifically, the output function $h:\mathbb{R}\hookrightarrow\Omega\times\mathbb{R}$ is defined by $x_t(\beta)\mapsto\big(\beta,x_t(\beta)\big)$, which gives a collection of aggregated measurements $Y_t=\{(\beta,y)\in\Omega\times\mathbb{R}: y_t=x_t(\beta)\text{ and } \b\in\Omega\}$ with $\b$ association with each measurement $y$. In this case, we may treat $Y_t=x_t\in L^2(\Omega)$ as a density function with respect to $\lambda$, which in turn defines a measure $\mu_t$ on $\Omega$ as $d\mu_t=x_td\lambda$. Applying the moment transform introduced in \eqref{eq:moment_integral}, we obtain $ m_k(t)=\int_{\Omega}\psi_kd\mu_t=\int_\Omega\psi_kx_td\lambda$.
Here we choose the monomial basis in $L^2(\Omega)$, i.e., $\psi_k=\beta^k$, which yields the moment dynamics,

\begin{align*}
    \frac{d}{dt}m_k(t)&=\frac{d}{dt}\int_0^1\beta^kx(t,\beta)d\beta=\int_0^1\beta^k\frac{d}{dt}x(t,\beta)d\beta \nonumber\\
    &=\int_0^1\beta^k\Big(\beta x(t,\beta)+\sum_{i=1}^p\beta^{i-1}u_i(t)\Big)d\beta\nonumber\\
    &=\int_0^1\beta^{k+1} x(t,\beta)d\beta+\sum_{i=1}^pu_i(t)\int_0^1\beta^{k+i-1}d\beta\nonumber\\
    &=m_{k+1}(t)+\sum_{i=1}^p\frac{1}{k+i}u_i(t). 
\end{align*}
This gives the moment system defined on the moment space $\cM$ of the form,
\begin{align}
    \frac{d}{dt}m(t)=Lm(t)+Hu(t), \label{eq:moment_linear_ensemble}
\end{align}
where $L:\mathcal{M}\rightarrow\mathcal{M}$ is the left-shift operator, given by $\big(m_0(t),m_1(t),\dots\big)'\mapsto\big(m_1(t),m_2(t),\dots\big)'$; $H:\mathbb{R}^p\rightarrow\mathcal{M}$ is a Hankel matrix with the $(k,i)$-entry equal to $\frac{1}{k+i}$ for $k=0,1,\dots$ and $i=1,\dots,p$; and $u(t)=\big(u_1(t),\dots,u_p(t)\big)'$. Consequently, the order-$q$ truncated moment system obeys 
\begin{align}
    \label{eq:truncated_moment_linear_ensemble}
    \frac{d}{dt}\widehat m^q(t)=\widehat L_q\widehat m^q(t)+\widehat H_qu(t) 
\end{align}
with
\begin{align*}
    \widehat L_q= \left[\begin{array}{cccc} 0 & 1 & & \\ & 0 & \ddots &  \\ & & \ddots & 1 \\ & & & 0 \end{array}\right],\widehat H_q=\left[\begin{array}{cccc} 1 & \frac{1}{2} &\cdots & \frac{1}{p} \\  \frac{1}{2} & \frac{1}{3} &\cdots & \frac{1}{p+1} \\ \vdots & \vdots & \ddots & \vdots \\ \frac{1}{q+1} & \frac{1}{q+2} & \cdots & \frac{1}{q+p} \end{array}\right].
\end{align*}
On the other end of the spectrum, we can compute the OT trajectory $\mu^*_t$ from $\mu_0$ and $\mu_1$ for $t\in [0,1]$ by using the DI formula in \eqref{eq:DI}. It follows that 
\begin{align*}
    m_k^*(t) &= \int_\Omega \big[(1-t)\beta+tX_1^{-1}\circ X_0(\beta)\big]^kd\mu_0(\beta)\\
    &=\int_\Omega \big[(1-t)\beta+tX_1^{-1}(X_0(\beta))\big]^kx_0(\beta)d\b,
\end{align*}
where $X_0$ and $X_1$ are the cumulative distribution functions of $\mu_0$ and $\mu_1$, respectively. 
Then, the moment dynamics describing this OT satisfy 
\begin{align*}
    &\frac{d}{dt}m_k^*(t)=\int_\Omega\frac{d}{dt}\big[(1-t)\beta+tX_1^{-1}(X_0(\beta))\big]^kx_0(\beta)d\beta\\
    &=\begin{cases}
        0,\qquad\qquad\qquad\qquad\qquad\qquad\qquad\qquad\qquad\, \text{if } k=0,\\
        \int_\Omega\big[X_1^{-1}(X_0(\beta))-\beta\big] \\ \quad\ \ \cdot\big[(1-t)\beta+tX_1^{-1}(X_0(\beta))\big]^{k-1}x_0(\beta)d\beta,\ \text{if }k\neq0.
    \end{cases}
\end{align*}

Following the formulation in \eqref{eq:OT_tracking_truncated}, the ensemble control input that transports the output measure $\mu_t$ from $\mu_0$ to $\mu_1$ can be found by solving the following OT-tracking control problem involving the order-$q$ truncation moment system, 
\begin{align}
    \min_{u:[0,1]\rightarrow\mathbb{R}^p}&\  \int_0^1 \|\widehat m^q(t)-P_qm^*(t)\|^2 dt\nonumber\\
    {\rm s.t.}\quad & \frac{d}{dt}\widehat m^q(t)=\widehat L_q\widehat m^q(t)+\widehat H_qu(t), \label{eq:OT_tracking_LS}
\end{align}
where $\|\cdot\|$ denotes the Euclidean norm on $\mathbb{R}^q$, and the truncated OT moment trajectory $P_qm^*(t)$ serves as the reference trajectory to be tracked. It can be shown that exact tracking is possible when ${\rm rank}(\widehat H_q)=q$, i.e., $p\geq q$ \cite{Brockett_1965}. This implies that when the number of control inputs is no less than the order of moment truncation, the tracking error can be made zero. The optimal tracking control can then be obtained, in feedback form, as the minimum-norm solution at each time $t\in [0,1]$, given by $u^*(t)=\widehat H_q'(\widehat H_q\widehat H_q')^{-1}\big(\frac{d}{dt}P_qm^*(t)-\widehat L_q\widehat m^q(t)\big)$. 

On the other hand, when $p<q$, exact tracking of the reference trajectory $P_qm^*(t)$ becomes infeasible. To strike a balance between the distributional control and optimal tracking tasks, we consider the fixed-endpoint control problem with a cost functional that includes a trade-off between tracking error and control energy as follow:
\begin{align}
   \min_{u:[0,1]\rightarrow\mathbb{R}^p} &\int_0^1\Big[ \|\widehat m^q(t)-P_qm^*(t)\|^2+u'(t)Ru(t)\Big]dt\nonumber\\
    {\rm s.t.} \ \ \ & \frac{d}{dt}\widehat m^q(t)=\widehat L_q\widehat m^q(t)+\widehat H_q u(t),\label{eq:tracking_control}\\
    & \widehat m^q(0)=P_qm^*(0),\ \ \widehat m^q(1)=P_qm^*(1),\nonumber
\end{align}
where $R\in\mathbb{R}^{p\times p}$ is a positive definite regulator. This optimal control problem can be solved using Pontryagin's maximum principle. The optimal control $u^*(t)=-\frac{1}{2}R^{-1}\widehat H_q'\lambda(t)$ is characterized by the costate $\lambda(t)$, which is a  solution to the following two-point boundary value problem \cite{Liberzon2011},
\begin{align*}
    \frac{d}{dt}\left[\begin{array}{c} \widehat m^q(t) \\ \lambda(t) \end{array}\right]  =\left[\begin{array}{cc} L & -\frac{1}{2}HR^{-1}H' \\ -2I & -L' \end{array}\right]&\left[\begin{array}{c} \widehat m^q(t) \\ \lambda(t) \end{array}\right]
     +2\left[\begin{array}{c} 0 \\ P_qm^*(t) \end{array}\right]
\end{align*}
with $\widehat m^q(0)=P_qm^*(0)$ and $\widehat m^q(1)=P_qm^*(1)$.

\begin{figure}[h]
\centering
\begin{subfigure}[b]{0.4\textwidth}
    \includegraphics[width=\linewidth]{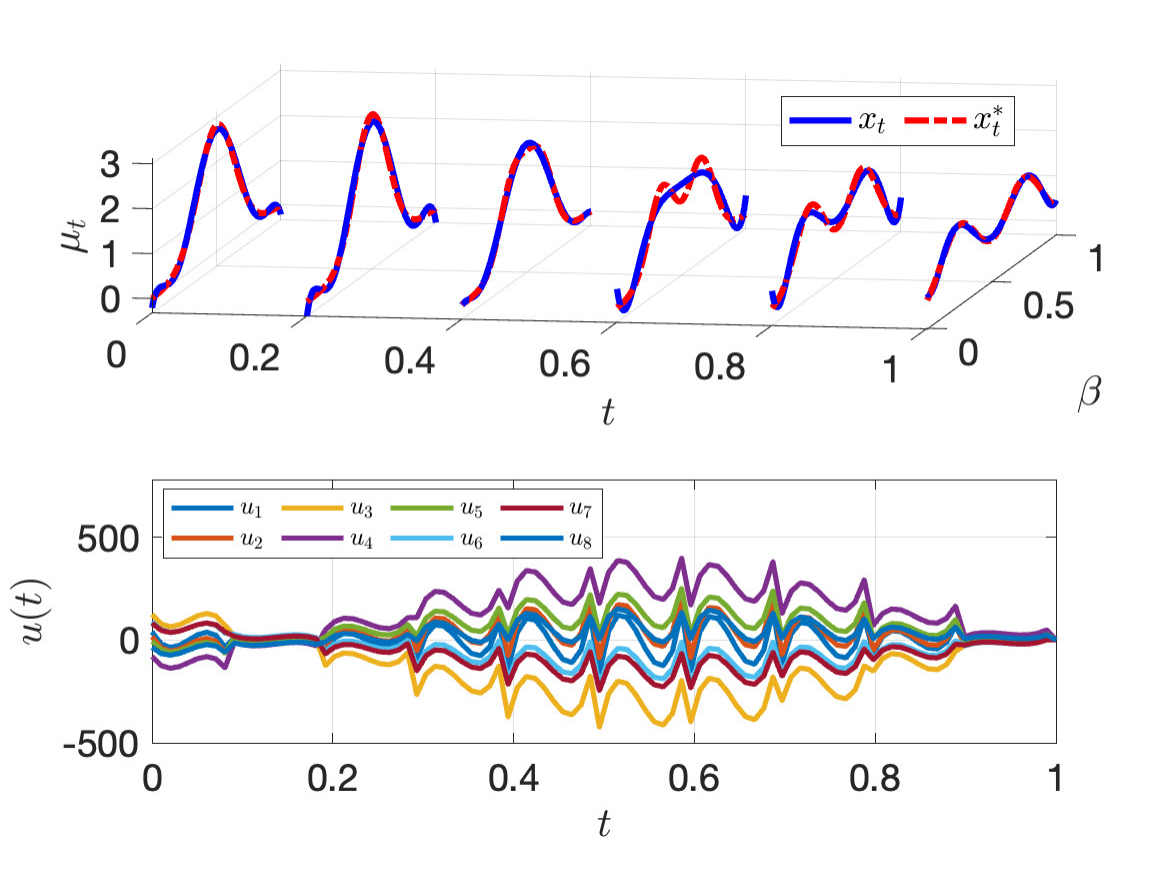}
    \caption{}
    \label{fig:linear_final}
\end{subfigure}
\begin{subfigure}[b]{0.4\textwidth}
    \includegraphics[width=\linewidth]{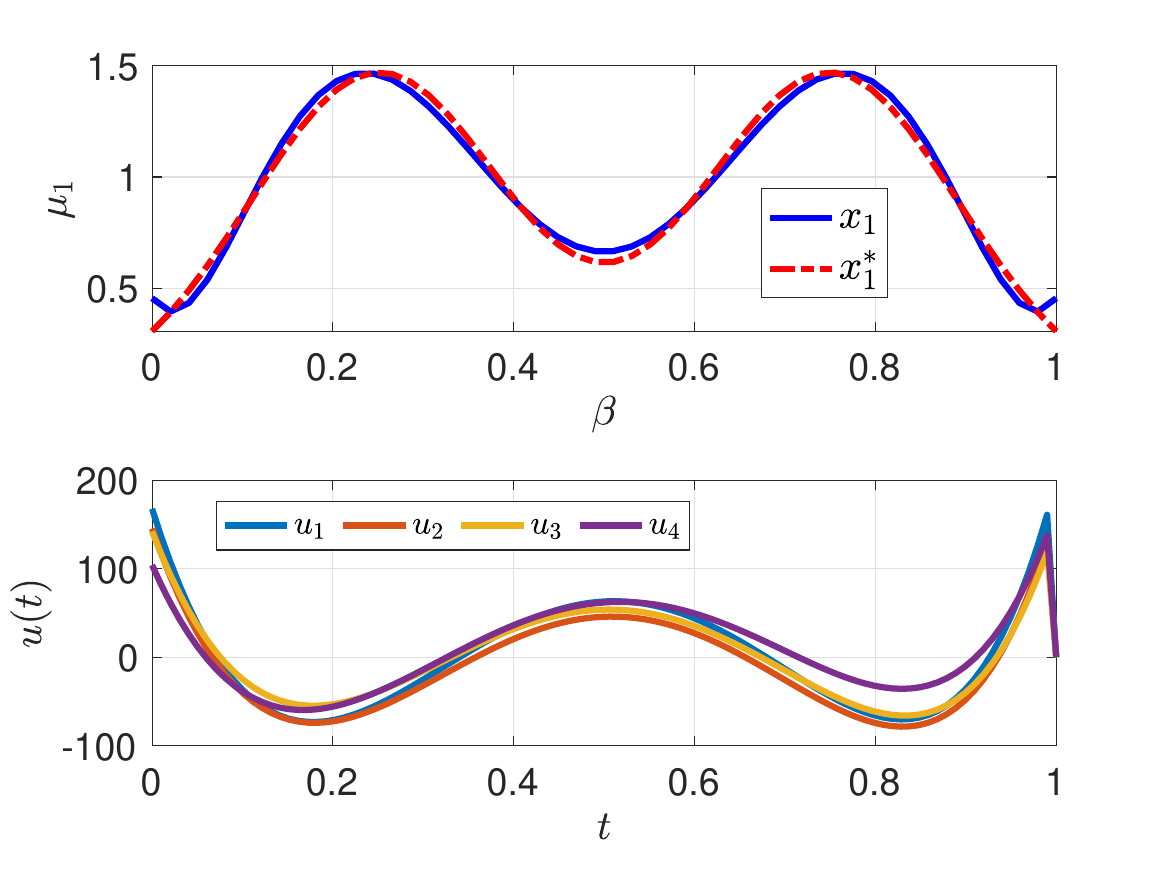}
        \caption{}
    \label{fig:linear_final_not_tracked}
\end{subfigure}
\caption{\small Functional control of the linear ensemble system $\Sigma_2$ in \eqref{eq:ECS_linear}. 
(a) plots the system trajectory $x_t$ (blue solid curves) and OT trajectory $x_t^*$ (red dashed curves) for $t=0,0.2,0.4,0.6,0.8,1$ (top panel) as well as the optimal tracking controls (bottom panel) for $q=8$ (moment truncation order) and $p=8$ (number of control inputs).
(b) shows the final ($x_1$, blue solid curve) and the desired ($x_1^*$, red dashed curve) state (top panel), and the calculated control inputs (bottom panel) for $q=8$ and $p=4$.
}
\label{fig:linear}
\end{figure}

To put this example into a concrete illustration, we consider steering the output measures of the ensemble $\Sigma_2$ from a truncated Gaussian density, $x_0^*(\b)=\psi\big(\frac{\b-a_0}{\sigma_0}\big)/\sigma_0\big[\Psi\big(\frac{1-a_0}{\sigma_0}\big)-\Psi\big(\frac{-a_0}{\sigma_0}\big)\big]$,
to a truncated Gaussian mixture density function, $x_1^*(\b)=\frac{c_{11}\psi(\a_{11})}{\sigma_{11}\big(\Psi(w_{11})-\Psi(z_{11})\big)}+\frac{c_{12}\psi(\a_{12})}{\sigma_{12}\big(\Psi(w_{12})-\Psi(z_{12})\big)}$.
Here, $\a_{ij}=\frac{\beta-a_{ij}}{\sigma_{ij}}$, $w_{ij}=\frac{1-a_{ij}}{\sigma_{ij}}$, and $z_{ij}=\frac{-a_{ij}}{\sigma_{ij}}$ for $i,j=1,2$; 
and $\psi$ and $\Psi$ are the probability density function and the cumulative distribution function of the standard normal distribution, respectively. We picked $a_0=0.5$, $a_{11}=0.25$, $a_{12}=0.75$, $\s_{0}=\s_{11}=\s_{12}=1/\sqrt{50}$, and $c_{11}=c_{12}=0.5$. The simulation results are shown in Fig. \ref{fig:linear}. In the case $p=q=8$, where the truncated OT trajectory can be perfectly tracked, we solved the OT-tracking control problem in \eqref{eq:OT_tracking_LS}. Fig. \ref{fig:linear_final} plots 
the optimal inputs $u_i(t)$, $t\in [0,1]$, for $i=1,\ldots,8$, and the controlled and OT trajectories $x_t$ and $x^*_t$, respectively, at sampled time instances, $t=0,0.2,0.4,0.6,0.8,1$. 
Similarly, in the case of $p=4$ and $q=8$, where the truncated OT trajectory cannot be perfectly tracked, Fig. \ref{fig:linear_final_not_tracked} shows the final and target density functions $x_1$ and $x^*_1$, respectively, and the optimal inputs obtained by solving the optimal control problem in \eqref{eq:tracking_control}.

Case II \emph{(Aggregated measurements without labels):} 
In numerous applications, recording the label of each individual system in an ensemble is impractical. In this case, the output function does not account for the system label $\b$ and is expressed as $h:\mathbb{R}\rightarrow\mathbb{R}$, mapping $x_t(\beta)$ to $h(x_t(\beta))$. To fix the idea, we also choose $h$ as the identity function, which gives the aggregated measurements $Y_t=\{y\in\mathbb{R}:y=x_t(\beta)\text{ and }\beta\in\Omega\}$. 
Using the the same set of basis functions $\psi_k(y)=y^k$, we obtain the $k^{\rm th}$ moment of the output measure $\mu_t$ induced by $Y_t$,
$$m_k(t)=\int_{\mathbb{R}} y^kd\mu_t(y)=\int_{\mathbb{R}}y^kd(x_t)_\#\lambda(y)=\int_0^1 x_t^k(\beta)d\beta$$
for $k\in\mathbb{N}$. Here, we consider driving $\S_2$ from the truncated Gaussian distribution ($\mu_0$) to the truncated Gaussian mixture distribution ($\mu_1$), with their probability density functions $x_0^*(\b)$ and $x_1^*(\b)$, respectively, as defined in Case I.
We solved the OT-tracking control problem in \eqref{eq:OT_tracking_truncated} for $p=q=8$, and the simulation results are shown in Fig. \ref{fig:linear_unlabeled}. Specifically, Fig. \ref{fig:unlabeled_flow} displays the empirical distribution approximation of $\mu_t$ for $t=0,0.2,0.4,0.6,0.8,1$ by using 1000 randomly selected individual systems in the ensemble, following the derived optimal control inputs $u_i$, $i=1,\ldots,8$ plotted in Figure \ref{fig:unlabeled_control}.

\begin{figure}[h]
\centering
\begin{subfigure}[b]{0.4\textwidth}
    \includegraphics[width=\linewidth]{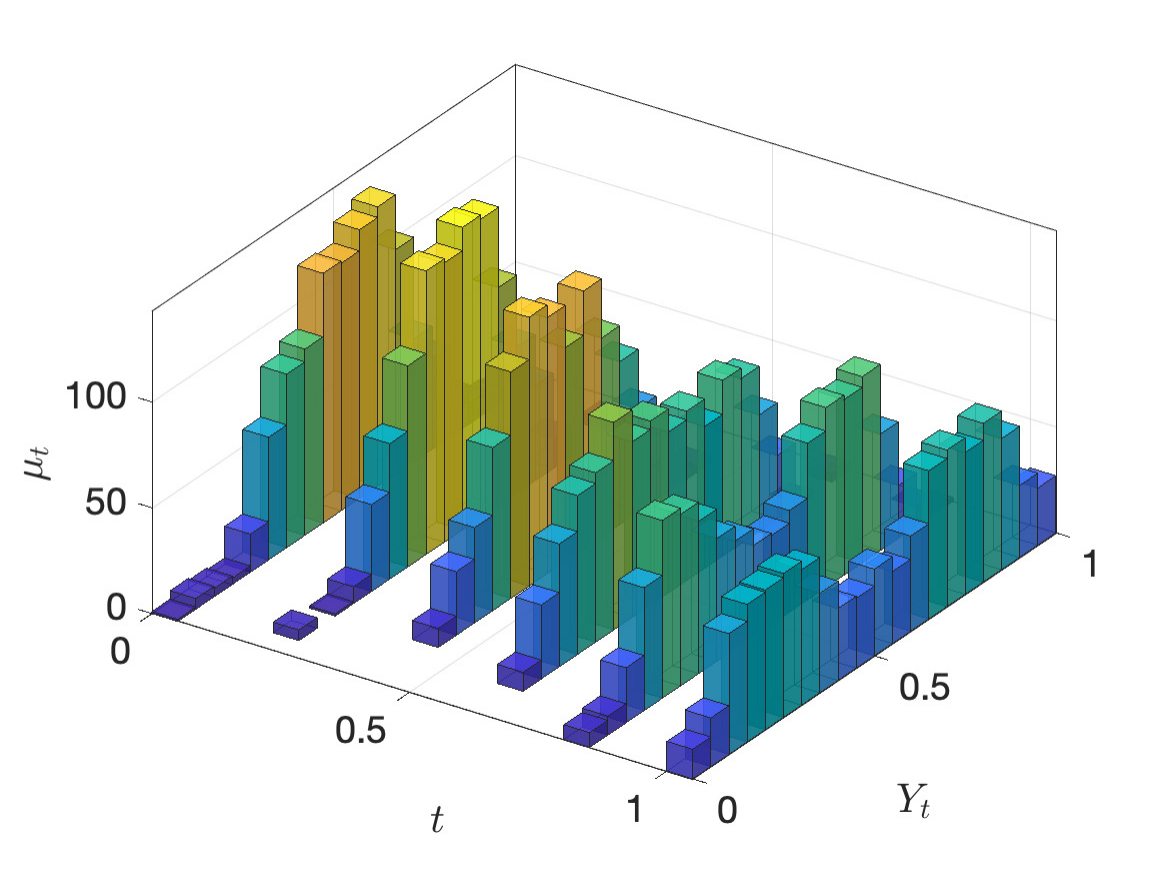}
    \caption{}
    \label{fig:unlabeled_flow}
\end{subfigure}
\begin{subfigure}[b]{0.4\textwidth}
    \includegraphics[width=\linewidth]{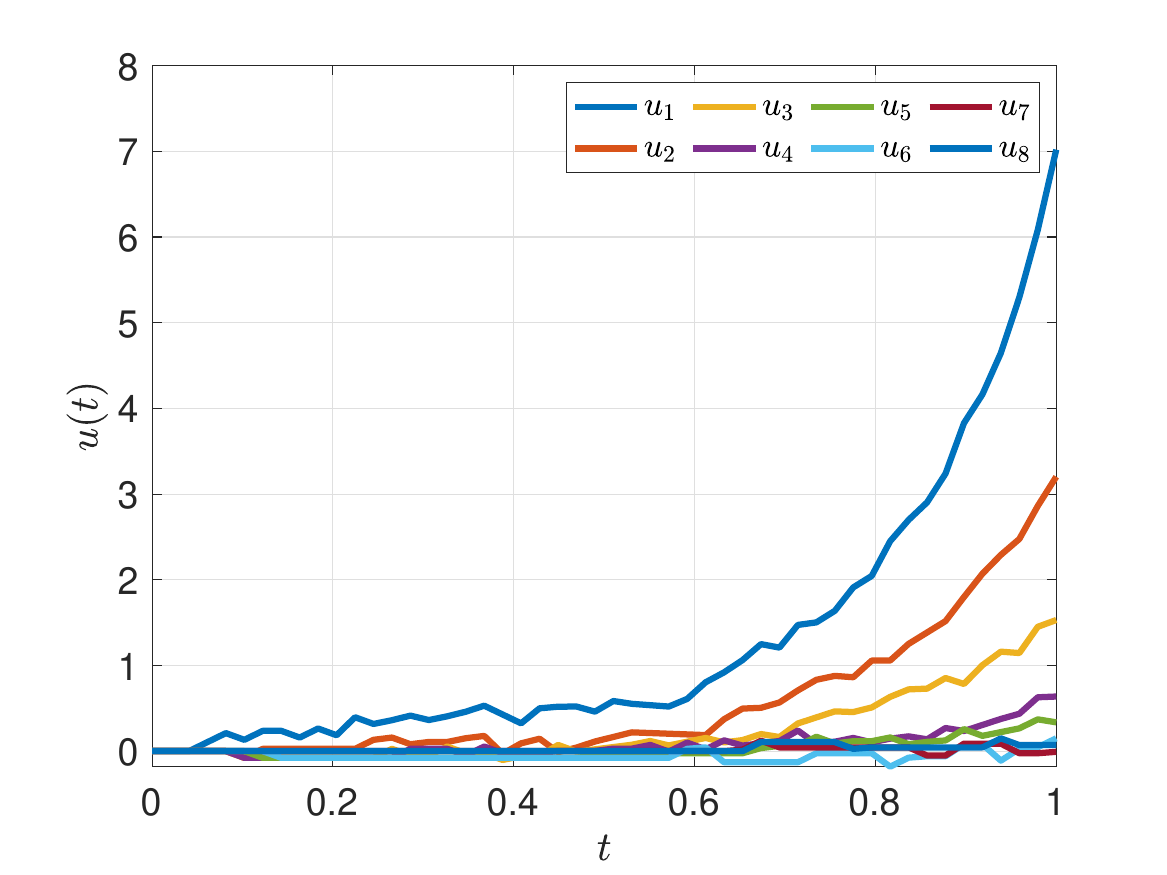}
    \caption{}
    \label{fig:unlabeled_control}
\end{subfigure}
\caption{\small Distributional control of the linear ensemble system $\Sigma_2$ in \eqref{eq:ECS_linear}. (a) shows the empirical distribution approximation of the controlled output measure $\mu_t$ at $t=0,0.2,0.4,0.6,0.8,1$ by using 1000 randomly selected individual systems in the ensemble, and (b) plots the control inputs for $q=8$ (moment truncation order) and $p=8$ (number of control inputs).}
\label{fig:linear_unlabeled}
\end{figure} 
\end{ex}

In the following, we present a distributional control problem frequently encountered in synchronization engineering to reinforce the applicability of the proposed method.

\begin{ex}[Synchronization engineering for rhythmic networks] 
\rm Synchronization engineering is concerned with coordinating a network of oscillatory systems to operate in unison \cite{Li_NatureComm16,Kiss2007}. In this example, we consider the problem of controlling synchronization in a large network of Kuramoto oscillators \cite{Strogatz2000}. Specifically, we consider a continuum of Kuramoto oscillators distributed on the circle $\mathbb{S}^1$ controlled by a common 
input $u(t)$, given by 
\begin{align}
    \frac{d}{dt}\theta(t,\w) &= \w+Kr(t)\sin(\psi(t)-\theta(t,\w))+u(t)\sin\theta(t,\w), \nonumber\\
    Y_t &= \t_t(\Omega). \label{eq:kuramoto}
\end{align}
Here, $\theta_t(\cdot)\doteq\theta(t,\cdot):\Omega\rightarrow\mathbb{S}^1$ denotes the phase of the oscillators, with their natural oscillation frequencies $\w$ distributed over a compact interval $\Omega\subset\mathbb{R}$; $K$ is the coupling strength; and $r(t)$ and $\psi(t)$ are the mean-field quantities determined by
$$r(t)e^{i\psi(t)}=\int_{\mathbb{S}^1}e^{i\theta}d\mu_t(\theta)=\int_{\Omega}e^{i\theta_t(\omega)}d\lambda(\w),$$ 
where $\lambda$ is the Haar (probability) measure on $\mathbb{S}^1$. In this case, the output measure $\mu_t=(\theta_t)_\#\lambda$ is a probability measure defined on $\mathbb{S}^1$. Engineering synchronization is equivalent to steering $\mu_t$ from an initial distribution $\mu_0$ 
to a point mass $\mu_1$, e.g., a $\delta$-distribution for complete synchrony, on $\mathbb{S}^1$ \cite{Li_NatureComm16}.

To tackle this distributional control problem, we define the moments using 
the Fourier basis $\psi_k(\theta)=e^{ik\theta}$ as these oscillators are periodic. 
This gives $m_k(t)=\<\psi_k,\mu_t\>=\int_{\mathbb{S}^1}e^{-ik\theta}d\mu_t(\theta)$ for $k\in\mathbb{N}$. 
In the simulation, we chose $\Omega=[-1,1]$ and the initial measure $\mu_0=\lambda$, which is the uniform distribution on $\mathbb{S}^1$. We then numerically solved the truncated OT-tracking problem in \eqref{eq:OT_tracking_truncated} with the truncation order $q=10$, and the simulation results are shown in Fig. \ref{fig:kuramoto}. Fig. \ref{fig:kuramoto_bar} displays the histograms of the empirical initial and final distributions generated using 1000 randomly selected oscillators in the Kuramoto ensemble. Fig. \ref{fig:kuramoto_distribution} plots the derived control input (top panel) and the respective density functions $\hat f_0$ and $\hat f_1$(with respect to $\lambda$) of the initial and final output measures, represented using order-10 truncated moment sequences given by
$ \hat f_{j}(\theta)=\frac{1}{2\pi}\big[\widehat m_0^{10}(j)+\sum_{k=1}^{10}\big((\widehat m_k^{10}(j))^\dagger e^{-ik\theta}+\widehat m_k^{10}(j)e^{ik\theta}\big)\big]$,
where $j=0,1$ and $(\widehat m_k^{10})^\dagger(j)$ denotes the complex conjugate of $\widehat m_k^{10}(j)$. The results demonstrate the ability to control the collective behavior of large-scale dynamic networks through the application of developed distributional control techniques.

\begin{figure}[h]
\centering
\begin{subfigure}[b]{0.4\textwidth}
    \includegraphics[width=\linewidth]{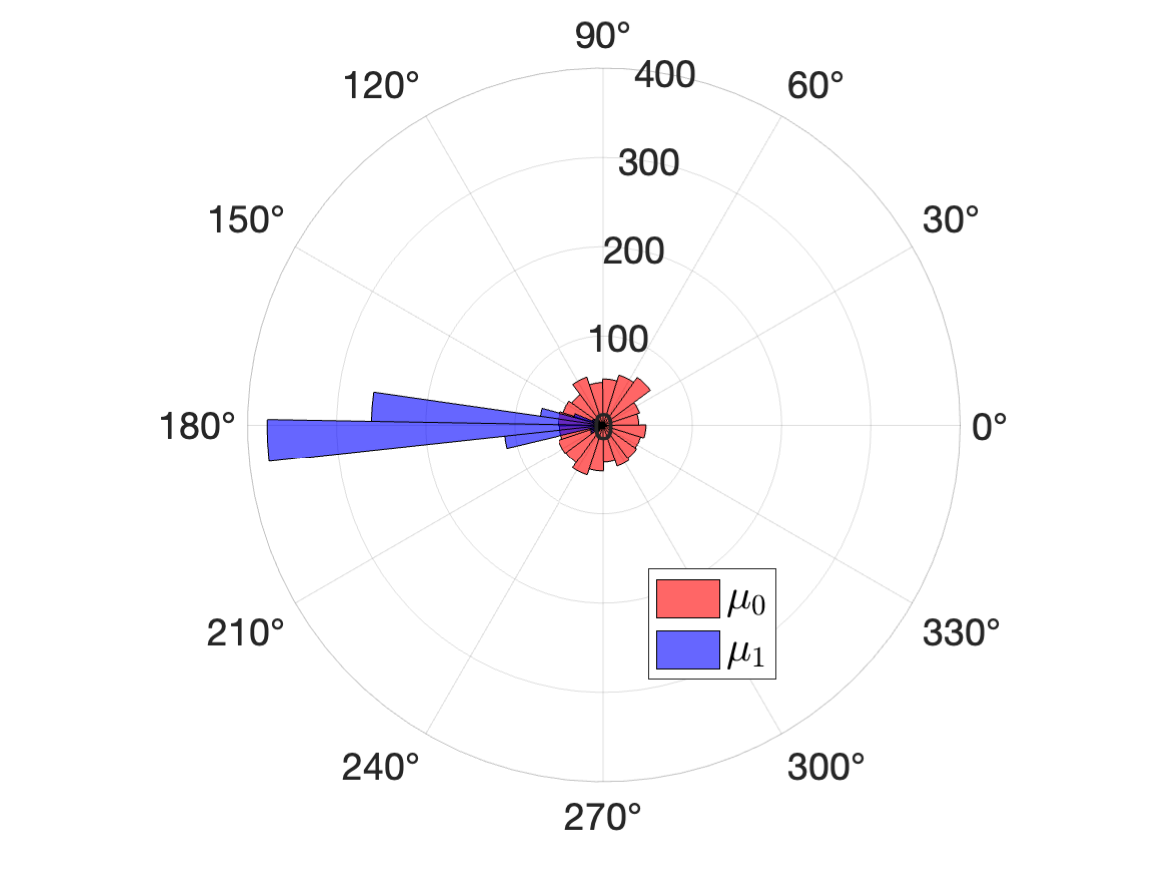}
    \caption{}
    \label{fig:kuramoto_bar}
\end{subfigure}
\begin{subfigure}[b]{0.4\textwidth}
    \includegraphics[width=\linewidth]{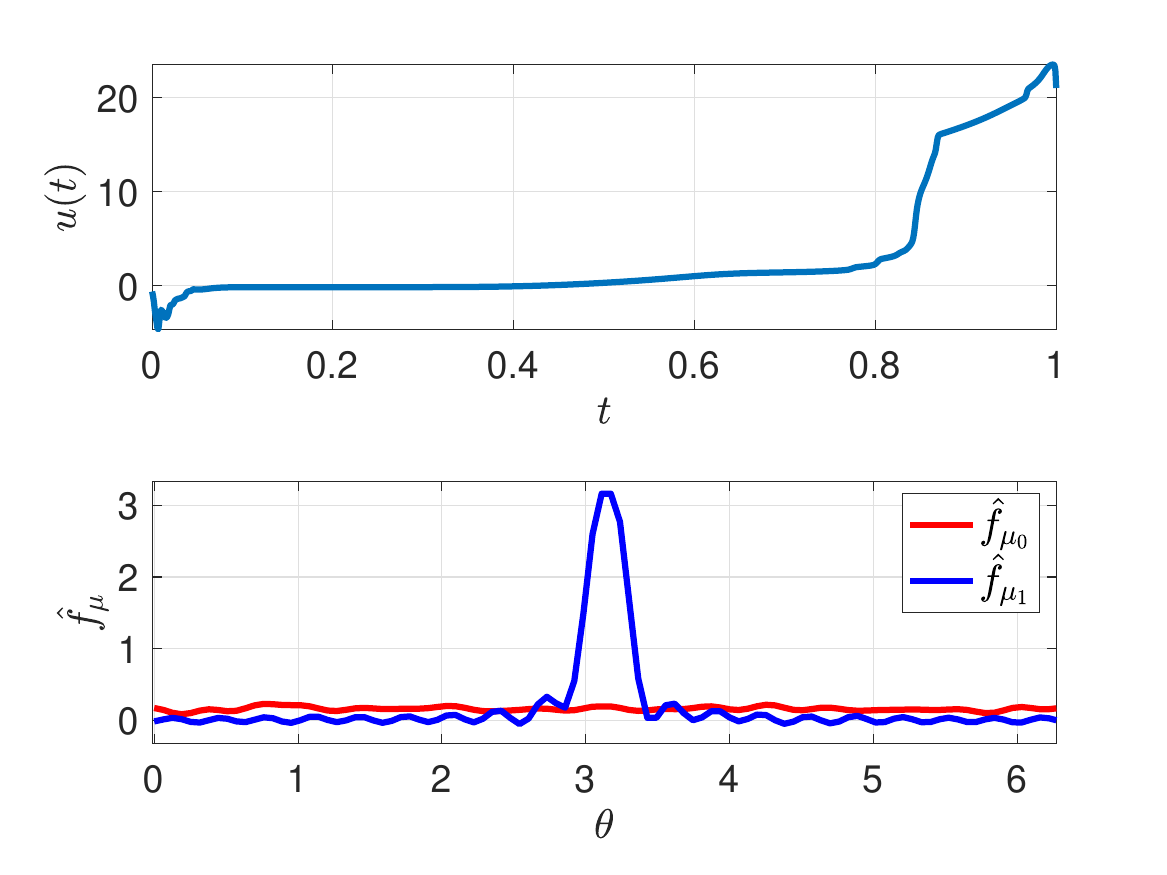}
    \caption{}
    \label{fig:kuramoto_distribution}
\end{subfigure}
\caption{\small Distributional control for synchronization of the Kuramoto oscillator network in \eqref{eq:kuramoto}. The initial and target distributions are the uniform distribution and a point mass on $\mathbb{S}^1$, respectively, with the moment truncation order chosen to be $q=10$. (a) shows the empirical distribution approximation of the initial (red) and final (blue) output measures generated by 1000 oscillators in the network. (b) plots the control input (top panel) and the approximated probability density functions of the initial and final distributions using the order-10 truncated moment sequences (bottom panel).}
\label{fig:kuramoto}
\end{figure}    
\end{ex}

\section{Conclusion}
\label{sec:conclusion}

In this paper, we introduce and formulate distributional control problems for ensemble systems. Shifting away from the conventional focus on controlling the states of an ensemble system, this new paradigm centers on controlling the time-varying output distributions induced by aggregated measurements of the ensemble system. We have developed a dynamic moment kernelization approach that enables systems-theoretic analysis and control design for distributional control problems. This emerging perspective on ensemble systems paves the way for further exploration and advancement in ensemble systems theory. In particular, it enables a purely data-driven paradigm of ensemble control, grounded in utilizing ensemble moment sequences that can be directly computed from available measurement data.

\appendices

\section*{Appendices}
\subsection{Duality between ${\it C_{\rm 0}(N)}$ and $\mathcal{P}({\it N})$}
\label{appd:measure}

Topologically, $N$ is a locally compact and separable Hausdorff space. As a result, every probability measure on $N$ is regular, and hence a Radon measure \cite{folland_13_real}. By the Riesz–Markov–Kakutani representation theorem, the map $\mu\mapsto I_\mu$, defined by $I_\mu(f)=\int_Nfd\mu$ is an isometric isomorphism from $M(N)$, the space of finite Radon measures on $N$, to $C_0^*(N)$, the dual of the space of real-valued functions on $N$ vanishing at infinity. As a result, the convergence $\mu_{t_n}\rightarrow\mu_t$ on $M(M)$ if and only if $\int_Nfd\mu_{t_n}\rightarrow\int_Nfd\mu_t$ for all $f\in C_0(N)$ generates the weak-$\star$ (vague) topology on $M(N)$, and hence on $\mathcal{P}(N)$. This topology on $\mathcal{P}(N)$ is generally weaker than the weak topology on $\mathcal{P}(N)$,  generated by the convergence $\mu_{t_n}\rightarrow\mu_t$ if and only if  $\int_Nfd\mu_{t_n}\rightarrow\int_N fd\mu_t$ for $f\in C_b(N)$, the space of bounded real-valued continuous functions on $N$. 
However, if the limit is a probability measure, then weak and weak-$\star$ convergence of the sequence coincides.

    \label{eq:moment_transform_1}

\begin{thm}
\label{thm:weak_convergence}
    Given a measure $\mu_t\in\mathcal{P}(N)$, then a sequence $\mu_{t_i}$ in $\mathcal{P}(N)$ satisfies $\mu_{t_i}\rightarrow\mu_t$ weakly as $i\rightarrow\infty$ if and only if $m_k(t_n)\rightarrow m_k(t)$ for all $k\in\mathbb{N}$.
\end{thm}
\begin{proof}
    We first assume that $\mu_{t_n}\rightarrow\mu_t$ weakly, then $\int_Nfd\mu_{t_n}\rightarrow\int_N fd\mu_t$ holds for all $f\in C_b(N)$. Because $C_0(N)\subseteq C_b(N)$, we obtain $m_k(t_n)\rightarrow m_k(t)$ for each $k\in\mathbb{N}$ by taking $f=\psi_k$. Conversely, we assume $m_k(t_n)\rightarrow m_k(t)$ for all $k\in\mathbb{N}$. Because $\{\psi_k\}_{k\in\mathbb{N}}$ is a basis for $C_0(N)$, for any $f\in C_0(N)$, there is a sequence $a_k$ in $\mathbb{R}$ such that $f=\sum_{k=0}^\infty a_k\psi_k$. By the dominant convergence theorem \cite{folland_13_real}, the integrability of $f$ implies $\sum_{k=0}^\infty a_km_k(t_n)=\sum_{k=0}^\infty\int_N\psi_kd\mu_{t_n}=\int_N \sum_{k=0}^\infty a_k\psi_kd\mu_{t_n}=\int_N fd\mu_{t_n}$ for all $n\in\mathbb{N}$, and the same result holds for $\mu_t$ and $m(t)$ as well. The application of dominant convergence theorem again further shows $\sum_{k=0}^\infty a_km_k(t_n)\rightarrow\sum_{k=0}^\infty a_km_k(t)$ as $n\rightarrow\infty$, yielding $\int_N fd\mu_{t_n}\rightarrow\int_Nfd\mu_t$ for any $f\in C_0(N)$, that is, the \emph{vague convergence} of $\mu_{t_n}$ to $\mu_t$. To show $\mu_{t_n}$ also converges to $\mu_t$ in the weak topology, we pick an arbitrary $f\in C_b(N)$ and denote its upper bound by $M$. For any $\varepsilon>0$, the regularity of $\mu$ indicates the existence a compact set $K\subset N$ such that $\mu(N\backslash K)=1-\mu(K)<\varepsilon$. Therefore, we can pick a continuous function $g$ supported on $K$ such that $0\leq g\leq1$ and $0<1-\int_Ngd\mu<\varepsilon/4M$. The vague convergence of $\mu_{t_n}$ to $\mu_t$ further implies that there is $n_0\in\mathbb{N}$ such that $0<1-\int_Ngd\mu_{t_n}<\varepsilon/3M$ for all $n>n_0$. Now, we have $gf\in C_0(N)$ so that $\int_Ngfd\mu_{t_n}\rightarrow\int_Ngfd\mu_t$ by using the vague convergence of $\mu_{t_n}$ to $\mu_t$, which particularly implies that the above $n_0$ can be chosen large enough to guarantee $\big|\int_Ngfd\mu_n-\int_Ngfd\mu\big|<\varepsilon/3$ for $n>n_0$. We then obtain the estimate $\big|\int_Nfd\mu_{t_n}-\int_Nfd\mu_t\big|\leq\big|\int_Nfd\mu_{t_n}-\int_Ngfd\mu_{t_n}\big|+\big|\int_Ngfd\mu_{t_n}-\int_Ngfd\mu_t\big|+\big|\int_Ngfd\mu_{t}-\int_Nfd\mu_t\big|\leq M\int_N(1-g)d\mu_{t_n}+\varepsilon/3+M\int_N(1-g)d\mu_{t}<\varepsilon$ for all $n>n_0$, concluding the weak convergence of $\mu_{t_n}$ to $\mu_t$.  
\end{proof}

\subsection{Continuity of the Output Measure Flow}
\label{appd:continuity}

 \begin{prop}
     \label{thm:distribution_dynamics}
     The output measure trajectory, i.e., the map $[0,T]\rightarrow\mathcal{P}(N)$ given by $t\mapsto\mu_t$, is continuous with respect to the weak topology on $\mathcal{P}(N)$, i.e., $\mu_t$ has a continuous representation.
 \end{prop}
 \begin{proof}
     The proof is based on first showing that $\mu_t:[0,T]\to\cP(N)$ is a continuous function with respect to the weak-$\star$ topology on $\cP(N)$, i.e., regarding $\mu_t$ as a Schwartz distribution. Then, we extend the weak-$\star$ continuity to weak continuity by showing that the family $\{\mu_t\}_{t\in[0,T]}$ of measures on $\mathcal{P}(N)$ indexed by $t\in[0,T]$ is tight.   
    To begin with, we first notice that Proposition \ref{prop:distribution_dynamics} (in the main text) actually shows that, for any $\varphi\in C_c^\infty(N)$, the function $[0,T]\rightarrow\mathbb{R}$ given by $t\mapsto\<\varphi,\mu_t\>$ is in $W^{1,1}([0,T])$, the Sobolev space consisting of integrable functions on $[0,T]$ with finite first-order weak derivatives, equivalently, the space of absolutely continuous functions on $[0,1]$ \cite{folland_13_real}. Let $L\subseteq[0,T]$ be the intersection of the Lebesgue sets of the collection of functions $\{t\mapsto\<\varphi_k,\mu_t\>\}_{k\in\mathbb{N}}$, then we have $\lambda(L)=1$, where $\{\varphi_k\}_{k\in\mathbb{N}}$ is a basis of $C_{c}^\infty(N)$ and $\lambda$ is the Lebesgue on $[0,1]$. Additionally, from \eqref{eq:distribution_dynamics}, we obtain for any $t_1,t_2\in L$, 
        $|\<\varphi_k,\mu_{t_1}-\mu_{t_2}\>|=\big|\int_{t_1}^{t_2}\<\varphi_k,\mu_{t}\>dt\big|
        \leq\int_{t_1}^{t_2}\int_\O |\nabla\varphi_k(y_t(\beta))||\nabla h(x_t(\beta))||F(t,\beta,x_t(\beta),u(t))|d\lambda(\beta)dt
        \leq C_h\sup_{y\in N}|\nabla\varphi_k(y)|\cdot$ \\ $\int_{t_1}^{t_2}\int_\Omega|F(t,\beta,x_t(\beta),u(t))|d\lambda(\beta)dt$, where $C_h$ is the Lipschitz constant of $h$. 
    Since $F$ is integrable, this implies that $|\<\varphi_k,\mu_{t_1}-\mu_{t_2}\>|\to 0$ uniformly in $k$ as $t_2\to t_1$. Therefore, $\{\mu_t\}_{t\in L}$ is a uniformly continuous family of continuous linear functionals on $C_c^\infty(N)$, and thus admits a unique extension to a continuous function $[0,T]\rightarrow\mathcal{P}(N)$, also denoted by $\mu_t$, with respect to the weak-* topology. To extend the weak-$\star$ continuity to weak continuity, it suffices to show the tightness of the family $\{\mu_t\}_{t\in[0,T]}$, ensuring that $\mu_{t_n}\rightarrow \mu_t$ weakly whenever $t_n\rightarrow t$. 
    
    To this end, we pick a sequence $\{U_k\}_{k=1}^\infty$ of precompact open subsets of $N$ such that $\overline{U}_k\subset U_{k+1}$ and $N=\cup_{i=1}^\infty U_i$. The existence of such a sequence is guaranteed by the local compactness and $\sigma$-compactness of $N$, where the $\sigma$-compactness is a result of the local compactness and separability of $N$ as a metric space \cite{folland_13_real}. We then choose a family of smooth functions $\{\phi_k\}_{k\in\mathbb{N}}$ satisfying $\phi_k(y)=1$ if $y\in\overline{U}_k$, $\phi_k(y)=0$ if $y\in N\backslash U_{k+1}$, and $|\nabla\phi_k|\leq 2/C_h$ for all $y\in N$, that is, a \emph{partition of unity} subordinate to $\{U_k\}_{k\in\mathbb{N}}$, whose existence is a direct consequence of Urysohn's lemma following from that $N$ is a locally compact Hausdorff space \cite{folland_13_real}. Applying the estimate derived above to $\phi_k$ yields $|\<\phi_k,\mu_{t}-\mu_{s}\>|\leq2\int_0^T\int_{y_\tau^{-1}(U_{k+1}\backslash\overline{U}_k)}|F(\tau,\beta,x_\tau(\beta),u(\tau))|d\lambda(\beta)d\tau$, and hence we obtain $\sum_{k\in\mathbb{N}}|\<\varphi_k,\mu_{t}-\mu_{s}\>|\leq2\int_0^T\int_{\Omega}|F(\tau,\beta,x_\tau(\beta),u(t))|d\lambda(\beta)d\tau<\infty$ by the integrability of $F$ again. Now, for a fixed $s\in[0,T]$ and any $\varepsilon>0$, because $\mu_s$ is inner regular, there exist $k\in\mathbb{N}$ such that $\<\varphi_k,\mu_s\>>1-\varepsilon/2$ and $\<\varphi_k,\mu_t-\mu_s\><\varepsilon/2$. This implies $\mu_t(\overline{U}_k)\geq\<\varphi_k,\mu_t\>\geq 1-\varepsilon$ for all $t\in[0,T]$, and therefore $\{\mu_t\}_{t\in[0,T]}$ is tight, concluding the proof.
\end{proof}

\subsection{Proof of Corollary \ref{cor:moment_convergence}}
\label{appd:proof_cor_2}

    To show the convergence of the sequence $J_q(u_q)$, we extend it to a double sequence $\big(J_p(u_q)\big)_{p,q\in\mathbb{N}}$. Here, $J_p(u_q)$ is defined as the minimal cost for the order-$q$ truncated moment system to track the order-$p$ truncated OT trajectory $P_pm^*(t)$, i.e., $J_p(u_q)=\min_{u}\Big\{\int_0^1d_{\cM}\big(P_pm^*(t),P_p\widehat m^q(t)\big)dt :\frac{d}{dt}\widehat m^q(t)= \widehat{F}^q(t,\widehat m^q(t),u(t))\Big\}$,  
    where, of course, $P_p\widehat m^q(t)=\widehat m^q(t)$ for $p>q$. Now, $\big(J_q(u_q)\big)_{q\in\mathbb{N}}$ embeds into the double sequence $\big(J_p(u_q)\big)_{p,q\in\mathbb{N}}$ as the diagonal subsequence. Consequently, if both the iterated limits $\lim_{p\rightarrow\infty}\lim_{q\rightarrow\infty}J_p(u_q)$ and $\lim_{q\rightarrow\infty}\lim_{p\rightarrow\infty}J_p(u_q)$ exists and equal to $J^*$, then $J_q(u_q)$ necessarily converges to $J^*$ by the diagonal argument \cite{Rudin1976}.

    We first compute $\lim_{p\rightarrow\infty}\lim_{ q\rightarrow\infty}J_p(u_q)$. For each fixed $p\in\mathbb{N}$, we observe that $J_p(u_q)$ is essentially the minimum of the set $$\Big\{\int_0^1d_{\cM}\big(P_pm^*(t),P_p\widehat m^{q+1}(t)\big)dt :\frac{d}{dt}\widehat m^q(t)= \widehat{F}^{q+1}(t,\widehat m^{q+1}(t),u(t)),\, \widehat m^{q+1}_{q+1}(t)=0\Big\}$$, a subset of $\Big\{\int_0^1d_{\cM}\big(P_pm^*(t),P_p\widehat m^{q+1}(t)\big)dt :\frac{d}{dt}\widehat m^q(t)= \widehat{F}^{q+1}(t,\widehat m^{q+1}(t),u(t))\Big\}$ whose minimum is $J_p(u_{q+1})$. This implies $J_p(u_{q})\geq J_p(u_{q+1})\geq 0$ for all $q\in\mathbb{N}$, yielding the convergence of the sequence $\big(J_p(u_q)\big)_{q\in\mathbb{N}}$.
    We denote the limit by $J_p(u)$, which is necessarily the minimal cost for the entire moment system to track $P_p m^*(t)$. Moreover, because $J_{p+1}(u)=\min_u\int_0^1 d_{\mathcal{M}}\big(P_{p+1}m^*(t), P_{p+1}m(t)\big)dt=\min_u\Big(\int_0^1 d_{\mathcal{M}}\big(P_pm^*(t),P_pm(t)\big)dt+\int_0^1d_{\cM}\big(P_{p+1}m^*(t)-P_pm^*(t),P_{p+1} m(t)-P_p m(t)\big)dt\Big)\geq J_p(u)+\min_u\int_0^1d_{\cM}\big(P_{p+1}m^*(t)-P_pm^*(t),P_{p+1} m(t)-P_p m(t)\big)dt$, the sequence $\big(J_p(u)\big)_{p\in\mathbb{N}}$ is monotonically increasing and bounded by $J^*$,  which gives $J(u)=\lim_{p\rightarrow\infty}\lim_{q\rightarrow\infty}J_p(q)=\lim_{p\rightarrow\infty}J_p(u)\leq J^*$. However, because $J^*$ is the minimum cost for the entire moment system to track $m^*(t)$, driven by any control input, the track cost must no less than $J^*$, particulatly $J(u)\geq J^*$. Therefore, we obtain $\lim_{p\rightarrow\infty}\lim_{q\rightarrow\infty}J_p(q)=J^*$.  

    On the other hand, for each fixed $q\in\mathbb{N}$, $\big(J_p(u_q)\big)_{p\geq q}$ is also monotonically increasing, following from $J_{p+1}(u_q)-J_p(u_q)=\int_0^1d_{\cM}(P_{p+1}m^*(t),P_pm^*(t))dt\geq0$, and bounded by $J(u_q)=J_q(u_q)+\int_0^1d_{\cM}(m^*(t),P_qm^*(t))dt$ so that $\lim_{p\rightarrow\infty} J_p(u_q)=J(u_q)$.
    For the sequence $\big(J(u_q)\big)_{q\in\mathbb{N}}$, following the same argument for $J_p(u_q)\geq J_p(u_{q+1})$ as above, $\big(J(u_q)\big)_{q\in\mathbb{N}}$ can also be shown to be a decreasing sequence lower bounded by $J^*$, and hence necessarily converges to $J^*$. This yields $\lim_{q\rightarrow\infty}\lim_{p\rightarrow\infty}J_p(u_q)=J^*$.

    Because both of the iterated limits of the double sequence $\big(J_p(u_q)\big)_{p,q\in\mathbb{N}}$ converge to the same limit $J^*$, the diagonal subsequence satisfies $J_q(u_q)\rightarrow J^*$ as desired.


\bibliographystyle{IEEEtran}
\bibliography{TAC_EC_OT_Refs}

\begin{thebibliography}{10}
\providecommand{\url}[1]{#1}
\csname url@samestyle\endcsname
\providecommand{\newblock}{\relax}
\providecommand{\bibinfo}[2]{#2}
\providecommand{\BIBentrySTDinterwordspacing}{\spaceskip=0pt\relax}
\providecommand{\BIBentryALTinterwordstretchfactor}{4}
\providecommand{\BIBentryALTinterwordspacing}{\spaceskip=\fontdimen2\font plus
\BIBentryALTinterwordstretchfactor\fontdimen3\font minus
  \fontdimen4\font\relax}
\providecommand{\BIBforeignlanguage}[2]{{%
\expandafter\ifx\csname l@#1\endcsname\relax
\typeout{** WARNING: IEEEtran.bst: No hyphenation pattern has been}%
\typeout{** loaded for the language `#1'. Using the pattern for}%
\typeout{** the default language instead.}%
\else
\language=\csname l@#1\endcsname
\fi
#2}}
\providecommand{\BIBdecl}{\relax}
\BIBdecl

\bibitem{Glaser98}
S.~J. Glaser, T.~Schulte-{Herbr\"uggen}, M.~Sieveking, N.~C.~N. O.~Schedletzky,
  O.~W. {S{\o}rensen}, and C.~Griesinger, ``Unitary control in quantum
  ensembles, maximizing signal intensity in coherent spectroscopy,''
  \emph{Science}, vol. 280, pp. 421--424, 1998.

\bibitem{Li_PNAS11}
J.-S. Li, J.~Ruths, T.-Y. Yu, H.~Arthanari, and G.~Wagner, ``Optimal pulse
  design in quantum control: A unified computational method,''
  \emph{Proceedings of the National Academy of Sciences}, vol. 108, no.~5, pp.
  1879--1884, 2011.

\bibitem{ching_13_control}
S.~Ching and J.~Ritt, ``Control strategies for underactuated neural ensembles
  driven by optogenetic stimulation,'' \emph{Frontiers in neural circuits},
  vol.~7, p.~54, 2013.

\bibitem{Kiss2007}
I.~Z. Kiss, C.~G. Rusin, H.~Kori, and J.~L. Hudson, ``Engineering complex
  dynamical structures: Sequential patterns and desynchronization,''
  \emph{Science}, vol. 316, no. 5833, pp. 1886--1889, 2007.

\bibitem{Sun2020}
C.~Sun, M.~Shen, and J.~P. How, ``Scaling up multiagent reinforcement learning
  for robotic systems: Learn an adaptive sparse communication graph,'' in
  \emph{2020 IEEE/RSJ International Conference on Intelligent Robots and
  Systems (IROS)}, 2020, pp. 11\,755--11\,762.

\bibitem{Li_TAC09}
J.-S. Li and N.~Khaneja, ``Ensemble control of \text{Bloch} equations,''
  \emph{IEEE Transactions on Automatic Control}, vol.~54, no.~3, pp. 528--536,
  2009.

\bibitem{Belhadj15}
M.~Belhadj, J.~Salomon, and G.~Turinici, ``Ensemble controllability and
  discrimination of perturbed bilinear control systems on connected, simple,
  compact lie groups,'' \emph{European Journal of Control}, vol.~22, pp.
  23--29, 2015.

\bibitem{Chen_Automatica19}
X.~Chen, ``Controllability of continuum ensemble of formation systems over
  directed graphs,'' \emph{Automatica}, vol. 108, p. 108497, 2019.

\bibitem{Boscain_SICON18}
N.~Augier, U.~Boscain, and M.~Sigalotti, ``Adiabatic ensemble control of a
  continuum of quantum systems,'' \emph{SIAM Journal on Control and
  Optimization}, vol.~56, no.~6, pp. 4045--4068, 2018.

\bibitem{Boscain_JDE22}
R.~Robin, N.~Augier, U.~Boscain, and M.~Sigalotti, ``Ensemble qubit
  controllability with a single control via adiabatic and rotating wave
  approximations,'' \emph{Journal of Differential Equations}, vol. 318, pp.
  414--442, 2022.

\bibitem{Dirr2016}
G.~Dirr, U.~Helmke, and M.~Schönlein, ``Controlling mean and variance in
  ensembles of linear systems,'' vol.~49, no.~18, pp. 1018--1023, 2016, 10th
  IFAC Symposium on Nonlinear Control Systems NOLCOS 2016.

\bibitem{Dirr2021}
G.~Dirr and M.~Sch{\"o}nlein, ``Uniform and ${L}_q$-ensemble reachability of
  parameter-dependent linear systems,'' \emph{Journal of Differential
  Equations}, vol. 283, pp. 216--262, 2021.

\bibitem{Schonlein2021}
M.~Sch{\"o}nlein, ``Ensemble reachability of homogenous parameter-depedent
  systems,'' \emph{PAMM}, vol.~20, no.~1, p. e202000342, 2021.

\bibitem{Schonlein2022}
------, ``Feedback equivalence and uniform ensemble reachability,''
  \emph{Linear Algebra and its Applications}, vol. 646, pp. 175--194, 2022.

\bibitem{Zeng_TAC15}
S.~Zeng, S.~Waldherr, C.~Ebenbauer, and F.~Allg{\"o}wer, ``Ensemble
  observability of linear systems,'' \emph{IEEE Transactions on Automatic
  Control}, vol.~61, no.~6, pp. 1452--1465, 2015.

\bibitem{Chen_Automatica20_Bloch}
X.~Chen, ``Ensemble observability of bloch equations with unknown population
  density,'' \emph{Automatica}, vol. 119, p. 109057, 2020.

\bibitem{Kiss2002}
I.~Z. Kiss, Y.~Zhai, and J.~L. Hudson, ``Emerging coherence in a population of
  chemical oscillators,'' \emph{Science}, vol. 296, no. 5573, pp. 1676--1678,
  2002.

\bibitem{Li_PRA_2006}
J.-S. Li and N.~Khaneja, ``Control of inhomogeneous quantum ensembles,''
  \emph{Phys. Rev. A}, vol.~73, p. 030302, Mar 2006.

\bibitem{Li_TAC11}
J.-S. Li, ``Ensemble control of finite-dimensional time-varying linear
  systems,'' \emph{IEEE Transactions on Automatic Control}, vol.~56, no.~2, pp.
  345--357, 2011.

\bibitem{Li_TAC13}
J.-S. Li, I.~Dasanayake, and J.~Ruths, ``Control and synchronization of neuron
  ensembles,'' \emph{IEEE Transactions on Automatic Control}, vol.~58, no.~8,
  pp. 1919--1930, 2013.

\bibitem{Helmke_SCL14}
U.~Helmke and M.~Sch{\"o}nlein, ``Uniform ensemble controllability for
  one-parameter families of time-invariant linear systems,'' \emph{Systems \&
  Control Letters}, vol.~71, pp. 69--77, 2014.

\bibitem{Zeng_SCL16}
S.~Zeng and F.~Allgöwer, ``A moment-based approach to ensemble controllability
  of linear systems,'' \emph{Systems \& Control Letters}, vol.~98, pp. 49--56,
  2016.

\bibitem{Chen_MCSS19}
X.~Chen, ``Structure theory for ensemble controllability, observability, and
  duality,'' \emph{Mathematics of Control, Signals, and Systems}, vol.~31,
  no.~2, p.~7, 2019.

\bibitem{Li_SICON21}
W.~Zhang and J.-S. Li, ``Ensemble control on {L}ie groups,'' \emph{SIAM Journal
  on Control and Optimization.}, vol.~59, no.~5, pp. 3805 -- 3827, 2021.

\bibitem{Zeng_TAC16}
S.~Zeng, H.~Ishii, and F.~Allg{\"o}wer, ``Sampled observability and state
  estimation of linear discrete ensembles,'' \emph{IEEE Transactions on
  Automatic Control}, vol.~62, no.~5, pp. 2406--2418, 2016.

\bibitem{Zeng_Automatica19}
S.~Zeng, ``Sample-based population observers,'' \emph{Automatica}, vol. 101,
  pp. 166--174, 2019.

\bibitem{Zeng_CDC16}
S.~Zeng and F.~Allg{\"o}wer, ``On the moment dynamics of discrete measures,''
  in \emph{2016 IEEE 55th Conference on Decision and Control}.\hskip 1em plus
  0.5em minus 0.4em\relax IEEE, 2016, pp. 4901--4906.

\bibitem{Miao_thesis}
W.~Miao, ``Algebraic, computational, and data-driven methods for
  control-theoretic analysis and learning of ensemble systems,'' 2021.

\bibitem{Krstic_TAC24}
V.~Alleaume and M.~Krstic, ``Ensembles of hyperbolic pdes: Stabilization by
  backstepping,'' \emph{IEEE Transactions on Automatic Control}, vol.~70,
  no.~2, pp. 905--920, 2025.

\bibitem{GRAPE}
N.~Khaneja, T.~Reiss, C.~Kehlet, T.~S.-Herbruggen, and S.~J. Glaser, ``Optimal
  control of coupled spin dynamics: design of \text{NMR} pulse sequences by
  gradient ascent algorithms,'' \emph{Journal of Magnetic Resonance}, vol. 172,
  pp. 296--305, 2005.

\bibitem{Li_JCP11}
J.~Ruths and J.-S. Li, ``A multidimensional pseudospectral method for optimal
  control of quantum ensembles,'' \emph{Journal of Chemical Physics}, vol. 134,
  p. 044128, 2011.

\bibitem{Dong_PRA14}
C.~Chen, D.~Dong, R.~Long, I.~R. Petersen, and H.~A. Rabitz, ``Sampling-based
  learning control of inhomogeneous quantum ensembles,'' \emph{Phys. Rev. A},
  vol.~89, no.~2, p. 023402, February 2014.

\bibitem{Li_CDC15_Fourier}
W.~Zhang and J.-S. Li, ``Uniform and selective excitations of spin ensembles
  with rf inhomogeneity,'' in \emph{2015 54th IEEE Conference on Decision and
  Control}.\hskip 1em plus 0.5em minus 0.4em\relax IEEE, 2015, pp. 5766--5771.

\bibitem{Gong_SICON16}
C.~Phelps, J.~O. Royset, and Q.~Gong, ``Optimal control of uncertain systems
  using sample average approximations,'' \emph{SIAM Journal on Control and
  Optimization}, vol.~54, no.~1, pp. 1--29, 2016.

\bibitem{Brent_JCP06}
B.~Pryor, ``Fourier decompositions and pulse sequence design algorithms for
  nuclear magnetic resonance in inhomogeneous fields,'' \emph{Journal of
  Chemical Physics}, vol. 125, p. 194111, 2006.

\bibitem{Li_ACC20}
W.~Miao and J.-S. Li, ``A geometric approach to linear ensemble control
  analysis and design,'' in \emph{2020 American Control Conference}, 2020, pp.
  4600--4605.

\bibitem{Li_ACC2024}
M.~Vu, B.~Singhal, J.-S. Li, and S.~Zeng, ``Data-driven moment-based control of
  linear ensemble systems,'' in \emph{2024 American Control Conference (ACC)},
  2024, pp. 5004--5009.

\bibitem{Zlotnik_ACC24}
A.~L.~P. de~Lima, A.~K. Harter, M.~J. Martin, and A.~Zlotnik, ``Optimal
  ensemble control of matter-wave splitting in bose-einstein condensates,'' in
  \emph{2024 American Control Conference (ACC)}, 2024, pp. 4196--4203.

\bibitem{Liberzon_SICON20}
D.~Liberzon and R.~W. Brockett, ``Spectral analysis of fokker--planck and
  related operators arising from linear stochastic differential equations,''
  \emph{SIAM Journal on Control and Optimization}, vol.~38, no.~5, pp.
  1453--1467, 2000.

\bibitem{Brockett2000}
R.~Brockett and N.~Khaneja, \emph{On the Stochastic Control of Quantum
  Ensembles}.\hskip 1em plus 0.5em minus 0.4em\relax Boston, MA: Springer US,
  2000, pp. 75--96.

\bibitem{Brockett2012}
R.~Brockett, \emph{Notes on the Control of the Liouville Equation}.\hskip 1em
  plus 0.5em minus 0.4em\relax Berlin, Heidelberg: Springer Berlin Heidelberg,
  2012, pp. 101--129.

\bibitem{Chen_TAC16_I}
Y.~Chen, T.~T. Georgiou, and M.~Pavon, ``Optimal steering of a linear
  stochastic system to a final probability distribution, part {I},'' \emph{IEEE
  Transactions on Automatic Control}, vol.~61, no.~5, pp. 1158--1169, 2016.

\bibitem{Chen_TAC16_II}
------, ``Optimal steering of a linear stochastic system to a final probability
  distribution, part {II},'' \emph{IEEE Transactions on Automatic Control},
  vol.~61, no.~5, pp. 1170--1180, 2016.

\bibitem{Chen_TAC16_OT}
------, ``Optimal transport over a linear dynamical system,'' \emph{IEEE
  Transactions on Automatic Control}, vol.~62, no.~5, pp. 2137--2152, 2016.

\bibitem{Becker_TRO12}
A.~Becker and T.~Bretl, ``Approximate steering of a unicycle under bounded
  model perturbation using ensemble control,'' \emph{IEEE Transactions on
  Robotics}, vol.~28, no.~3, pp. 580--591, 2012.

\bibitem{Narayanan_ACC19}
V.~Narayanan, J.~T. Ritt, J.-S. Li, and S.~Ching, ``A learning framework for
  controlling spiking neural networks,'' in \emph{2019 American Control
  Conference (ACC)}, 2019, pp. 211--216.

\bibitem{Li_ACC20_Learning}
Y.-C. Yu, V.~Narayanan, S.~Ching, and J.-S. Li, ``Learning to control neurons
  using aggregated measurements,'' in \emph{2020 American Control Conference},
  2020, pp. 4028--4033.

\bibitem{Dong_PRA22}
C.~Jiang, Y.~Pan, Z.-G. Wu, Q.~Gao, and D.~Dong, ``Robust optimization for
  quantum reinforcement learning control using partial observations,''
  \emph{Phys. Rev. A}, vol. 105, p. 062443, Jun 2022.

\bibitem{Li_NatureComm16}
A.~Zlotnik, R.~Nagao, I.~Z. Kiss, and J.-S. Li, ``Phase-selective entrainment
  of nonlinear oscillator ensembles,'' \emph{Nature Communications}, vol.~7, p.
  10788, 2016.

\bibitem{Bogachev07}
V.~I. Bogachev, \emph{Measure Theory}.\hskip 1em plus 0.5em minus 0.4em\relax
  Springer-Verlag Berlin Heidelberg, 2007.

\bibitem{Strogatz2000}
S.~H. Strogatz, ``From kuramoto to crawford: exploring the onset of
  synchronization in populations of coupled oscillators,'' \emph{Physica D:
  Nonlinear Phenomena}, vol. 143, no.~1, pp. 1--20, 2000.

\bibitem{Silver85}
M.~S. Silver, R.~I. Joseph, and D.~I. Hoult, ``Selective spin inversion in
  nuclear magnetic resonance and coherent optics through an exact solution of
  the bloch-riccati equation,'' \emph{Physical Review A}, vol.~31, no.~4, pp.
  2753--2755, 1985.

\bibitem{folland_13_real}
G.~B. Folland, \emph{Real analysis: modern techniques and their
  applications}.\hskip 1em plus 0.5em minus 0.4em\relax John Wiley \& Sons,
  2013.

\bibitem{villani2009}
C.~Villani, \emph{Optimal transport: old and new}.\hskip 1em plus 0.5em minus
  0.4em\relax Springer, 2009, vol. 338.

\bibitem{Billingsley95}
P.~Billingsley, \emph{Probability and Measure}, 3rd~ed., ser. Wiley Series in
  Probability and Statistics.\hskip 1em plus 0.5em minus 0.4em\relax Wiley,
  1995, vol. 245.

\bibitem{brockett_2000_stochastic}
R.~Brockett and N.~Khaneja, ``On the stochastic control of quantum ensembles,''
  in \emph{System theory}.\hskip 1em plus 0.5em minus 0.4em\relax Springer,
  2000, pp. 75--96.

\bibitem{hausdorff_23_momentprobleme}
F.~Hausdorff, ``Momentprobleme f{\"u}r ein endliches intervall.''
  \emph{Mathematische Zeitschrift}, vol.~16, no.~1, pp. 220--248, 1923.

\bibitem{McCann97}
\BIBentryALTinterwordspacing
R.~J. McCann, ``A convexity principle for interacting gases,'' \emph{Advances
  in Mathematics}, vol. 128, no.~1, pp. 153--179, 1997. [Online]. Available:
  \url{https://www.sciencedirect.com/science/article/pii/S0001870897916340}
\BIBentrySTDinterwordspacing

\bibitem{zeng_18_computation}
S.~Zeng, W.~Zhang, and J.-S. Li, ``On the computation of control inputs for
  linear ensembles,'' in \emph{2018 Annual American Control Conference
  (ACC)}.\hskip 1em plus 0.5em minus 0.4em\relax IEEE, 2018, pp. 6101--6107.

\bibitem{Li_CDC24}
Y.-H. Kuan, W.~Zhang, and J.-S. Li, ``Computational moment control of ensemble
  systems,'' in \emph{The 63rd IEEE Conference on Decision and Control}, Milan,
  Italy, December 2024.

\bibitem{Brockett_1965}
R.~Brockett and M.~Mesarovi{\'c}, ``The reproducibility of multivariable
  systems,'' \emph{Journal of Mathematical Analysis and Applications}, vol.~11,
  pp. 548--563, 1965.

\bibitem{Liberzon2011}
D.~Liberzon, \emph{Calculus of Variations and Optimal Control Theory: A Concise
  Introduction}.\hskip 1em plus 0.5em minus 0.4em\relax Princeton University
  Press, 2011.

\bibitem{Rudin1976}
W.~Rudin, \emph{Principles of Mathematical Analysis}.\hskip 1em plus 0.5em
  minus 0.4em\relax McGraw-Hill, 1976.

\end{thebibliography}

\end{document}